\documentclass[reqno]{amsart}

\usepackage{amsmath,amsfonts,amssymb,amsthm}
\usepackage{enumitem,physics,cite,color,setspace}
\usepackage{tikz}
\usepackage{subcaption}
\usepackage{hyperref}

\numberwithin{equation}{section}
\allowdisplaybreaks
\makeatletter
\def\subsection{\@startsection{subsection}{2}%
  \z@{.5\baselineskip}%
  {.5\baselineskip}%
  {\normalfont\centering\scshape}%
}
\makeatother

\newtheorem{theorem}{Theorem}[section]
\newtheorem{lemma}{Lemma}[section]
\newtheorem{cor}{Corollary}[theorem]
\newtheorem{example}{Example}[section]
\newtheorem{definition}{Definition}[section]
\newtheorem{remark}{Remark}[section]

\newcommand{\e}{\mathrm{e}}
\newcommand{\Srec}{\mathcal{S}}      
\newcommand{\Trec}{\mathcal{T}}      

\title[Branch points and natural boundaries of solutions of FDEs]{Branch points and natural boundaries of Poincaré functions and of solutions of nonlinear functional differential equations}
\author[Techheang Meng]{Techheang Meng$^*$}
\address{$^*$Department of Mathematics, China University of Mining and Technology-Beijing, Beijing 100083, P. R. China}
\email{$^*$mengtechheangbbb@yahoo.com}

\newcommand{\figPoincareOne}{
\begin{figure}[htbp]
     \centering
\hspace{-2.5cm}
     \begin{subfigure}[b]{0.3\textwidth}
         \centering

         \caption{First generation}
         \label{fig:poinc:1a}
     \end{subfigure}
     \hspace{2.5cm}
     \begin{subfigure}[b]{0.3\textwidth}
         \centering

         \caption{Second iteration}
         \label{fig:poinc:1b}
     \end{subfigure}
     \caption{Branch point pattern of (\ref{eq:poinc:2})}
     \label{fig:poinc:1}
\end{figure}
}

\newcommand{\figPoincareTwo}{
\begin{figure}[htbp]
     \centering
     \hspace{-2.5cm}
     \begin{subfigure}[b]{0.3\textwidth}
         \centering

\caption{First generation}
         \label{fig:poinc:2a}
     \end{subfigure}
     \hspace{2.5cm}
     \begin{subfigure}[b]{0.3\textwidth}
         \centering

\caption{Second generation}
         \label{fig:poinc:2b}
     \end{subfigure}
     \caption{Branch point patterns of (\ref{eq:poinc:6})}
     \label{fig:poinc:2}
\end{figure}
}

\newcommand{\figPoincareThree}{
\begin{figure}[htbp]
     \centering
     \hspace{-2.5cm}
     \begin{subfigure}[b]{0.3\textwidth}
         \centering

         \caption{First generation}
         \label{fig:poinc:3a}
     \end{subfigure}
     \hspace{2.5cm}
     \begin{subfigure}[b]{0.3\textwidth}
         \centering

         \caption{Second generation}
         \label{fig:poinc:3b}
     \end{subfigure}
     \caption{Branch point patterns of (\ref{eq:poinc:10})}
     \label{fig:poinc:3}
\end{figure}
}

\newcommand{\figPoincareFour}{
\begin{figure}[htbp]
    \centering

\begin{tikzpicture}[x=1.2pt,y=1.2pt,yscale=-1,xscale=1]

\draw    (272.54,71.98) -- (297.19,71.98) ;
\draw    (297.19,71.98) -- (335.29,71.98) ;
\draw    (335.29,71.98) -- (354.22,71.98) ;
\draw    (207.6,71.98) -- (234.45,71.98) ;
\draw [shift={(234.45,71.98)}, rotate = 0] [color={rgb, 255:red, 0; green, 0; blue, 0 }  ][fill={rgb, 255:red, 0; green, 0; blue, 0 }  ][line width=0.75]      (0, 0) circle [x radius= 2.01, y radius= 2.01]   ;
\draw    (234.45,71.98) -- (272.54,71.98) ;
\draw [shift={(272.54,71.98)}, rotate = 0] [color={rgb, 255:red, 0; green, 0; blue, 0 }  ][fill={rgb, 255:red, 0; green, 0; blue, 0 }  ][line width=0.75]      (0, 0) circle [x radius= 2.01, y radius= 2.01]   ;
\draw  [fill={rgb, 255:red, 208; green, 2; blue, 27 }  ,fill opacity=1 ] (332.95,71.98) .. controls (332.95,70.69) and (334,69.64) .. (335.29,69.64) .. controls (336.58,69.64) and (337.62,70.69) .. (337.62,71.98) .. controls (337.62,73.27) and (336.58,74.31) .. (335.29,74.31) .. controls (334,74.31) and (332.95,73.27) .. (332.95,71.98) -- cycle ;

\draw [color={rgb, 255:red, 245; green, 166; blue, 35 }  ,draw opacity=1 ]   (272.54,121.78) -- (297.19,121.78) ;
\draw [color={rgb, 255:red, 245; green, 166; blue, 35 }  ,draw opacity=1 ]   (297.19,121.78) -- (335.29,121.78) ;
\draw [color={rgb, 255:red, 245; green, 166; blue, 35 }  ,draw opacity=1 ]   (335.29,121.78) -- (354.22,121.78) ;
\draw [color={rgb, 255:red, 245; green, 166; blue, 35 }  ,draw opacity=1 ]   (207.6,121.78) -- (234.45,121.78) ;
\draw [color={rgb, 255:red, 245; green, 166; blue, 35 }  ,draw opacity=1 ]   (234.45,121.78) -- (272.54,121.78) ;
\draw  [color={rgb, 255:red, 0; green, 0; blue, 0 }  ,draw opacity=1 ][fill={rgb, 255:red, 208; green, 2; blue, 27 }  ,fill opacity=1 ] (332.95,121.78) .. controls (332.95,120.49) and (334,119.44) .. (335.29,119.44) .. controls (336.58,119.44) and (337.62,120.49) .. (337.62,121.78) .. controls (337.62,123.07) and (336.58,124.11) .. (335.29,124.11) .. controls (334,124.11) and (332.95,123.07) .. (332.95,121.78) -- cycle ;
\draw  [fill={rgb, 255:red, 74; green, 144; blue, 226 }  ,fill opacity=1 ] (232.11,121.78) .. controls (232.11,120.49) and (233.16,119.44) .. (234.45,119.44) .. controls (235.74,119.44) and (236.78,120.49) .. (236.78,121.78) .. controls (236.78,123.07) and (235.74,124.11) .. (234.45,124.11) .. controls (233.16,124.11) and (232.11,123.07) .. (232.11,121.78) -- cycle ;
\draw [color={rgb, 255:red, 144; green, 19; blue, 254 }  ,draw opacity=1 ]   (272.54,171.58) -- (297.19,171.58) ;
\draw [color={rgb, 255:red, 144; green, 19; blue, 254 }  ,draw opacity=1 ]   (297.19,171.58) -- (335.29,171.58) ;
\draw [color={rgb, 255:red, 144; green, 19; blue, 254 }  ,draw opacity=1 ]   (335.29,171.58) -- (354.22,171.58) ;
\draw [color={rgb, 255:red, 144; green, 19; blue, 254 }  ,draw opacity=1 ]   (207.6,171.58) -- (234.45,171.58) ;
\draw [color={rgb, 255:red, 144; green, 19; blue, 254 }  ,draw opacity=1 ]   (234.45,171.58) -- (272.54,171.58) ;
\draw  [color={rgb, 255:red, 0; green, 0; blue, 0 }  ,draw opacity=1 ][fill={rgb, 255:red, 208; green, 2; blue, 27 }  ,fill opacity=1 ] (332.95,171.58) .. controls (332.95,170.29) and (334,169.24) .. (335.29,169.24) .. controls (336.58,169.24) and (337.62,170.29) .. (337.62,171.58) .. controls (337.62,172.87) and (336.58,173.91) .. (335.29,173.91) .. controls (334,173.91) and (332.95,172.87) .. (332.95,171.58) -- cycle ;
\draw  [color={rgb, 255:red, 0; green, 0; blue, 0 }  ,draw opacity=1 ][fill={rgb, 255:red, 74; green, 144; blue, 226 }  ,fill opacity=1 ] (232.11,171.58) .. controls (232.11,170.29) and (233.16,169.24) .. (234.45,169.24) .. controls (235.74,169.24) and (236.78,170.29) .. (236.78,171.58) .. controls (236.78,172.87) and (235.74,173.91) .. (234.45,173.91) .. controls (233.16,173.91) and (232.11,172.87) .. (232.11,171.58) -- cycle ;

\draw (231.05,76.18) node [anchor=north west][inner sep=0.75pt]  [font=\scriptsize] [align=left] {$\displaystyle 0$};
\draw (255.49,76.98) node [anchor=north west][inner sep=0.75pt]  [font=\scriptsize,color={rgb, 255:red, 208; green, 2; blue, 27 }  ,opacity=1 ] [align=left] {$\displaystyle 4/729$};
\draw (318.24,74.98) node [anchor=north west][inner sep=0.75pt]  [font=\scriptsize,color={rgb, 255:red, 208; green, 2; blue, 27 }  ,opacity=1 ] [align=left] {$\displaystyle \frac{1}{9}( \times 2)$};
\draw (217.02,91.98) node [anchor=north west][inner sep=0.75pt]  [font=\scriptsize,color={rgb, 255:red, 74; green, 144; blue, 226 }  ,opacity=1 ] [align=left] {$\displaystyle \frac{1}{3}( \times 2)$};
\draw (321.89,156.98) node [anchor=north west][inner sep=0.75pt]  [font=\scriptsize,color={rgb, 255:red, 208; green, 2; blue, 27 }  ,opacity=1 ] [align=left] {$\displaystyle 4/9$};
\draw (318.24,124.78) node [anchor=north west][inner sep=0.75pt]  [font=\scriptsize,color={rgb, 255:red, 208; green, 2; blue, 27 }  ,opacity=1 ] [align=left] {$\displaystyle \frac{1}{9}( \times 2)$};
\draw (213.82,141.58) node [anchor=north west][inner sep=0.75pt]  [font=\scriptsize,color={rgb, 255:red, 74; green, 144; blue, 226 }  ,opacity=1 ] [align=left] {$\displaystyle \frac{1}{3}( \times 2)$};

\end{tikzpicture}

    \caption{Branch point pattern of (\ref{eq:poinc:14})}
    \label{fig:poinc:4}
\end{figure}
}

\newcommand{\figPoincareFive}{
\begin{figure}[htbp]
     \centering

\begin{tikzpicture}[x=1pt,y=1pt,yscale=-1,xscale=1]

\draw    (335.05,75.18) -- (359.7,75.18) ;
\draw [shift={(359.7,75.18)}, rotate = 0] [color={rgb, 255:red, 0; green, 0; blue, 0 }  ][fill={rgb, 255:red, 0; green, 0; blue, 0 }  ][line width=0.75]      (0, 0) circle [x radius= 2.01, y radius= 2.01]   ;
\draw    (409.06,75.18) -- (447.15,75.18) ;
\draw    (209.5,75.18) -- (236.35,75.18) ;
\draw [shift={(236.35,75.18)}, rotate = 0] [color={rgb, 255:red, 0; green, 0; blue, 0 }  ][fill={rgb, 255:red, 0; green, 0; blue, 0 }  ][line width=0.75]      (0, 0) circle [x radius= 2.01, y radius= 2.01]   ;
\draw [shift={(209.5,75.18)}, rotate = 0] [color={rgb, 255:red, 0; green, 0; blue, 0 }  ][fill={rgb, 255:red, 0; green, 0; blue, 0 }  ][line width=0.75]      (0, 0) circle [x radius= 2.01, y radius= 2.01]   ;
\draw    (236.35,75.18) -- (285.7,75.18) ;
\draw [shift={(285.7,75.18)}, rotate = 0] [color={rgb, 255:red, 0; green, 0; blue, 0 }  ][fill={rgb, 255:red, 0; green, 0; blue, 0 }  ][line width=0.75]      (0, 0) circle [x radius= 2.01, y radius= 2.01]   ;
\draw    (285.7,75.18) -- (335.05,75.18) ;
\draw    (359.7,75.18) -- (409.06,75.18) ;
\draw  [fill={rgb, 255:red, 208; green, 2; blue, 27 }  ,fill opacity=1 ] (409.06,75.18) .. controls (409.06,73.89) and (410.1,72.84) .. (411.39,72.84) .. controls (412.68,72.84) and (413.73,73.89) .. (413.73,75.18) .. controls (413.73,76.47) and (412.68,77.51) .. (411.39,77.51) .. controls (410.1,77.51) and (409.06,76.47) .. (409.06,75.18) -- cycle ;
\draw    (190.57,75.18) -- (209.5,75.18) ;

\draw [color={rgb, 255:red, 245; green, 166; blue, 35 }  ,draw opacity=1 ]   (335.05,117.51) -- (359.7,117.51) ;
\draw [color={rgb, 255:red, 245; green, 166; blue, 35 }  ,draw opacity=1 ]   (409.06,117.51) -- (447.15,117.51) ;
\draw [color={rgb, 255:red, 245; green, 166; blue, 35 }  ,draw opacity=1 ]   (209.5,117.51) -- (236.35,117.51) ;
\draw [color={rgb, 255:red, 245; green, 166; blue, 35 }  ,draw opacity=1 ]   (236.35,117.51) -- (285.7,117.51) ;
\draw [color={rgb, 255:red, 245; green, 166; blue, 35 }  ,draw opacity=1 ]   (285.7,117.51) -- (335.05,117.51) ;
\draw [color={rgb, 255:red, 245; green, 166; blue, 35 }  ,draw opacity=1 ]   (359.7,117.51) -- (409.06,117.51) ;
\draw  [color={rgb, 255:red, 0; green, 0; blue, 0 }  ,draw opacity=1 ][fill={rgb, 255:red, 208; green, 2; blue, 27 }  ,fill opacity=1 ] (409.06,117.51) .. controls (409.06,116.22) and (410.1,115.17) .. (411.39,115.17) .. controls (412.68,115.17) and (413.73,116.22) .. (413.73,117.51) .. controls (413.73,118.8) and (412.68,119.84) .. (411.39,119.84) .. controls (410.1,119.84) and (409.06,118.8) .. (409.06,117.51) -- cycle ;
\draw [color={rgb, 255:red, 245; green, 166; blue, 35 }  ,draw opacity=1 ]   (190.57,117.51) -- (209.5,117.51) ;
\draw  [color={rgb, 255:red, 0; green, 0; blue, 0 }  ,draw opacity=1 ][fill={rgb, 255:red, 74; green, 144; blue, 226 }  ,fill opacity=1 ] (357.37,117.51) .. controls (357.37,116.22) and (358.41,115.17) .. (359.7,115.17) .. controls (360.99,115.17) and (362.04,116.22) .. (362.04,117.51) .. controls (362.04,118.8) and (360.99,119.84) .. (359.7,119.84) .. controls (358.41,119.84) and (357.37,118.8) .. (357.37,117.51) -- cycle ;

\draw [color={rgb, 255:red, 144; green, 19; blue, 254 }  ,draw opacity=1 ]   (335.05,159.84) -- (359.7,159.84) ;
\draw [color={rgb, 255:red, 144; green, 19; blue, 254 }  ,draw opacity=1 ]   (409.06,159.84) -- (447.15,159.84) ;
\draw [color={rgb, 255:red, 144; green, 19; blue, 254 }  ,draw opacity=1 ]   (209.5,159.84) -- (236.35,159.84) ;
\draw [color={rgb, 255:red, 144; green, 19; blue, 254 }  ,draw opacity=1 ]   (236.35,159.84) -- (285.7,159.84) ;
\draw [color={rgb, 255:red, 144; green, 19; blue, 254 }  ,draw opacity=1 ]   (285.7,159.84) -- (335.05,159.84) ;
\draw [color={rgb, 255:red, 144; green, 19; blue, 254 }  ,draw opacity=1 ]   (359.7,159.84) -- (409.06,159.84) ;
\draw  [color={rgb, 255:red, 0; green, 0; blue, 0 }  ,draw opacity=1 ][fill={rgb, 255:red, 80; green, 227; blue, 194 }  ,fill opacity=1 ] (409.06,159.84) .. controls (409.06,158.55) and (410.1,157.5) .. (411.39,157.5) .. controls (412.68,157.5) and (413.73,158.55) .. (413.73,159.84) .. controls (413.73,161.13) and (412.68,162.17) .. (411.39,162.17) .. controls (410.1,162.17) and (409.06,161.13) .. (409.06,159.84) -- cycle ;
\draw [color={rgb, 255:red, 144; green, 19; blue, 254 }  ,draw opacity=1 ]   (190.57,159.84) -- (209.5,159.84) ;
\draw  [color={rgb, 255:red, 0; green, 0; blue, 0 }  ,draw opacity=1 ][fill={rgb, 255:red, 74; green, 144; blue, 226 }  ,fill opacity=1 ] (357.37,159.84) .. controls (357.37,158.55) and (358.41,157.5) .. (359.7,157.5) .. controls (360.99,157.5) and (362.04,158.55) .. (362.04,159.84) .. controls (362.04,161.13) and (360.99,162.17) .. (359.7,162.17) .. controls (358.41,162.17) and (357.37,161.13) .. (357.37,159.84) -- cycle ;

\draw [color={rgb, 255:red, 139; green, 87; blue, 42 }  ,draw opacity=1 ]   (335.05,202.18) -- (359.7,202.18) ;
\draw [shift={(359.7,202.18)}, rotate = 0] [color={rgb, 255:red, 139; green, 87; blue, 42 }  ,draw opacity=1 ][fill={rgb, 255:red, 139; green, 87; blue, 42 }  ,fill opacity=1 ][line width=0.75]      (0, 0) circle [x radius= 2.01, y radius= 2.01]   ;
\draw [color={rgb, 255:red, 139; green, 87; blue, 42 }  ,draw opacity=1 ]   (409.06,202.18) -- (447.15,202.18) ;
\draw [color={rgb, 255:red, 139; green, 87; blue, 42 }  ,draw opacity=1 ]   (209.5,202.18) -- (236.35,202.18) ;
\draw [color={rgb, 255:red, 139; green, 87; blue, 42 }  ,draw opacity=1 ]   (236.35,202.18) -- (285.7,202.18) ;
\draw [color={rgb, 255:red, 139; green, 87; blue, 42 }  ,draw opacity=1 ]   (285.7,202.18) -- (335.05,202.18) ;
\draw [color={rgb, 255:red, 139; green, 87; blue, 42 }  ,draw opacity=1 ]   (359.7,202.18) -- (409.06,202.18) ;
\draw  [color={rgb, 255:red, 0; green, 0; blue, 0 }  ,draw opacity=1 ][fill={rgb, 255:red, 80; green, 227; blue, 194 }  ,fill opacity=1 ] (409.06,202.18) .. controls (409.06,200.89) and (410.1,199.84) .. (411.39,199.84) .. controls (412.68,199.84) and (413.73,200.89) .. (413.73,202.18) .. controls (413.73,203.47) and (412.68,204.51) .. (411.39,204.51) .. controls (410.1,204.51) and (409.06,203.47) .. (409.06,202.18) -- cycle ;
\draw [color={rgb, 255:red, 139; green, 87; blue, 42 }  ,draw opacity=1 ]   (190.57,202.18) -- (209.5,202.18) ;

\draw (281.45,78.58) node [anchor=north west][inner sep=0.75pt]  [font=\scriptsize] [align=left] {$\displaystyle 0$};
\draw (192.57,78.18) node [anchor=north west][inner sep=0.75pt]  [font=\scriptsize,color={rgb, 255:red, 208; green, 2; blue, 27 }  ,opacity=1 ] [align=left] {$\displaystyle -1/36$};
\draw (341.02,119.98) node [anchor=north west][inner sep=0.75pt]  [font=\scriptsize,color={rgb, 255:red, 74; green, 144; blue, 226 }  ,opacity=1 ] [align=left] {$\displaystyle \frac{5}{12}( \times 2)$};
\draw (395.84,77.58) node [anchor=north west][inner sep=0.75pt]  [font=\scriptsize,color={rgb, 255:red, 208; green, 2; blue, 27 }  ,opacity=1 ] [align=left] {$\displaystyle \frac{1}{3}( \times 2)$};
\draw (211.42,59.98) node [anchor=north west][inner sep=0.75pt]  [font=\scriptsize,color={rgb, 255:red, 74; green, 144; blue, 226 }  ,opacity=1 ] [align=left] {$\displaystyle -575/20736$};
\draw (349.02,60.78) node [anchor=north west][inner sep=0.75pt]  [font=\scriptsize,color={rgb, 255:red, 74; green, 144; blue, 226 }  ,opacity=1 ] [align=left] {$\displaystyle 1/6$};
\draw (395.04,120.78) node [anchor=north west][inner sep=0.75pt]  [font=\scriptsize,color={rgb, 255:red, 208; green, 2; blue, 27 }  ,opacity=1 ] [align=left] {$\displaystyle \frac{1}{3}( \times 2)$};
\draw (397.82,162.38) node [anchor=north west][inner sep=0.75pt]  [font=\scriptsize,color={rgb, 255:red, 80; green, 227; blue, 194 }  ,opacity=1 ] [align=left] {$\displaystyle \frac{1}{2}( \times 2)$};
\draw (398.62,204.78) node [anchor=north west][inner sep=0.75pt]  [font=\scriptsize,color={rgb, 255:red, 80; green, 227; blue, 194 }  ,opacity=1 ] [align=left] {$\displaystyle \frac{1}{2}( \times 2)$};
\draw (312.38,208.18) node [anchor=north west][inner sep=0.75pt]  [font=\scriptsize,color={rgb, 255:red, 74; green, 144; blue, 226 }  ,opacity=1 ] [align=left] {$\displaystyle \left( 5+2\sqrt{3}\right) /12$};
\draw (337.05,162.84) node [anchor=north west][inner sep=0.75pt]  [font=\scriptsize,color={rgb, 255:red, 74; green, 144; blue, 226 }  ,opacity=1 ] [align=left] {$\displaystyle \frac{5}{12}( \times 2)$};

\end{tikzpicture}

\caption{Branch point pattern of (\ref{eq:poinc:16})}
    \label{fig:poinc:5}
\end{figure}
}

\begin{document}

\begin{abstract}
In this paper, we study the singularity structure of the analytic continuation of solutions to the equation $g(\lambda w)=R(g(w))$, where $|\lambda|<1$ and $R$ is a rational function. We show that the analytic solution \(g\) near zero is branched, and that its analytic continuation yields an infinite-sheeted Riemann surface with various patterns of branch points. In one specific pattern, some branch points are delayed: they appear only when the analytic continuation reaches a particular Riemann sheet. Using this functional-equation framework, we extend the study to first-order nonlinear functional differential equations. We show that the continuation of certain special solutions near an attracting fixed point yields either single-valued or multivalued functions, with the unit circle as a natural boundary on every sheet.
\end{abstract}

\maketitle

\section{Introduction}
\label{sec:intro}

In complex dynamics, the behaviour of a rational function $R$ near some fixed points, namely the attracting and the repelling ones, can be studied using the Poincar\'e function\footnote{The inverse of this function is also known as the Schr\"oder or Koenigs function \cite[Ch.3, §4]{steinmetz2011rational}.} \cite{koenigs:84,poincare:90,steinmetz2011rational,carleson2013complex}. A fixed point $t_0$ of $R$ is called attracting if $0<|R'(t_0)|<1$, superattracting if $R'(t_0)=0$, repelling if $|R'(t_0)|>1$, rationally neutral if $R'(t_0)$ is a root of unity, and irrationally neutral otherwise. After a M\"obius change of variables, we may assume that $t_0=0$. If $0$ is neither neutral nor superattracting, the Poincar\'e function $g$ is the solution of
\begin{equation}
    \label{eq:intro:1}
    g(\lambda w)=R(g(w)),\quad g(0)=0,\quad g'(0)=1,
\end{equation}
where, after a further M\"obius transformation, $R$ can be normalised to the form
\begin{equation}
    \label{eq:intro:2}
    R(t)=\frac{\lambda t+a_2t^2+\cdots+a_pt^p}{1+b_1t+\cdots+b_qt^q}.
\end{equation}
Here, $|\lambda|\neq0$ and $|\lambda|\neq1$ since $0$ is neither neutral nor superattracting. The role of the function $g$ implies that the dynamics of $R$ near $0$ can be understood similarly as a dynamics of a linear function $t\mapsto \lambda t$ at $0$, see \cite[Ch.3, §3]{steinmetz2011rational}. This existence of such solutions was shown by Koenigs \cite{koenigs:84}. When $0$ is a repelling fixed point, Poincar\'e proved that $g$ can be analytically continued to a meromorphic function on the whole complex plane \cite{poincare:90}. In this case, $y(z)=g(\lambda^z)$ is a meromorphic solution of a difference equation
\begin{equation}
\label{eq:intro:3}
	y(z+1)=R(y(z)),
\end{equation}
normalised by $y(z)\lambda^{-z}\to \alpha$ as $\Re z\to-\infty$. Meromorphic solutions of (\ref{eq:intro:3}) have been completely described in a series of works \cite{p2kimura1971,p2Simo,p2Yana,p2BakerLiverpool,p2Azarina}. The study of the equation (\ref{eq:intro:3}), and is also motivated from different perspectives in discrete integrable systems, see \cite{veselov1991integrable,halburd2006}.

When $0$ is an attracting fixed point, the Poincar\'e function still exists, but its analytic continuation has branch points. In a joint work with Halburd, Korhonen, and Liu \cite{halburd2025extensionanalyticsolutionsfirstorder}, we described the singularity structure completely for the quadratic family
\begin{equation}
    \label{eq:intro:4}
R(t)=\lambda t+t^2,\quad  0<\lambda<1.
\end{equation}
This description is closely related to the work of Mahler \cite{mahler1976class,mahler1981special} on the equation
\begin{equation}
    \label{eq:intro:5}
    Y(w^2)=Y(w)^2+c,\quad Y(w)=\frac{1}{w}+\sum_{k\geq0}a_kw^k,
\end{equation}
whose solution always has the unit circle as a natural boundary. Indeed, $y(z)=Y(\exp{-2^{z}})-\lambda/2$ satisfies (\ref{eq:intro:3}) with $R$ as in (\ref{eq:intro:4}) and $c=\lambda/2-\lambda^2/4$; this is a different solution from the normalized one above, but it solves the same difference equation. As observed in \cite[p.149]{meng2026application}, the function $Y$ can be expressed in terms of a B\"ottcher function (see \cite[Ch.3 §2]{steinmetz2011rational}) of a superattracting fixed point at $\infty$.

In \cite{marshall2002natural,marshall2004natural}, Marshall, van Brunt and Wake studied the linear pantograph-type equation (see \cite{iserles1993generalized,derfel1997pantograph} for its applications)
\begin{equation}
    \label{eq:intro:6}
    y''(z)+ay'(z)+by(z)=cy(h(z)),
\end{equation}
near a (super)attracting fixed point $z_0$ of an entire function $h$. They proved that an analytic solution at $z_0$ can be continued to the immediate basin of attraction $A_{z_0}^0$, the connected component containing $z_0$ of the set of points whose orbits under $h$ converging to $z_0$, and that $\partial A_{z_0}^0$ is a natural boundary. Their proof combines the iteration theory of entire functions with the growth estimates of entire functions.

Motivated by these works of Marshall, van Brunt and Wake as well as the idea of Mahler on finding a solution with a pole, we study the nonlinear FDE
\begin{equation}
\label{eq:intro:7}
y'+ay^2+by(z^2)=0,
\end{equation}
where $a,b$ are complex constants and $y$ has a pole at $z=0$. Although (\ref{eq:intro:7}) contains only one nonlinear term, its solutions need not be single-valued on $A_0^0\setminus\{0\}$. We will show that every transcendental solution has the natural boundary $\partial A_0^0=C(0,1)$ on each Riemann sheet, and that both single-valued and multivalued solutions occur, depending on the radius of convergence of the analytic part of the Laurent series at $0$ (Theorem \ref{thm:fde:1}). 

As a multiplicative analogue of (\ref{eq:intro:1}), we also study the $q$-difference equation:
\begin{equation}
\label{eq:intro:8}
y'+ay+by^2(qz)=0,\quad q\neq0
\end{equation}
where $y$ has a pole at $z=0$. Here the analytic continuation is more rigid and depends on $q$ (Theorem \ref{thm:qde:2}). When $|q|<1$, the solutions are meromorphic. When $|q|=1$ and $q$ is not a root of unity, they have a natural boundary. When $|q|>1$, the continuation must have branch points, whose structure can be described as in \cite{halburd2025extensionanalyticsolutionsfirstorder}. When $q$ is a root of unity, (\ref{eq:intro:8}) reduces to a system of ODEs which in general is not integrable\footnote{in a sense of not having enough first integrals}, and we will discuss its meromorphic solutions using arguments similar to those in \cite{halburd2015all,meng2024meromorphicsolutions3dlotkavolterra}.

This work started from Mahler's equation (\ref{eq:intro:5}) and the results of Marshall, van Brunt and Wake on (\ref{eq:intro:6}). Our goal was to see if we can describe a singularity structure of a solution with a pole at an attracting fixed point of any simple nonliner FDE. The Riccati-type equations (\ref{eq:intro:7}) and (\ref{eq:intro:8}) are among the simplest containing both ideas of works. The quadratic case of the Poincar\'e function was then treated in the joint work \cite{halburd2025extensionanalyticsolutionsfirstorder} with Halburd, Korhonen and Liu. After that paper, the author developed a familiarity and thus, the language of \emph{types}, \emph{stems}, \emph{sheets} and \emph{generations}, which describes multivalued continuations. Since this language is very clear for Poincar\'e functions, we present it first, in Section~\ref{sec:poinc}, and use it again in Sections~\ref{sec:fde} and~\ref{sec:qde}. While working on the higher degree Poincar\'e functions, we found that the first branch point of a given type can be \emph{delayed}: it may first appear on a sheet of a later generation, rather than on the germ $g_0$ of $g$ at $0$.

In Section~\ref{sec:poinc}, we extend the results of \cite{halburd2025extensionanalyticsolutionsfirstorder} to rational functions $R$ of higher degree and study the patterns formed by the branch points. There are two distinct patterns, depending on whether the stem of a given type lies on $g_0$. If it does, as for the quadratic polynomial (\ref{eq:intro:4}), then $g_0$ has infinitely many branch points of that type, located at $w_0\lambda^{1-k}$, $k\in\mathbb{N}$. If it does not, then each sheet carrying the stem has exactly one branch point of that type, and some sheets have none (Lemma~\ref{lem:poinc:3}). A complete classification of the possible patterns does not seem feasible; instead, the examples in Section~\ref{subsec:poinc:2} show the variety that occurs. In particular, the stem can first appear on the $k$-th generation for any prescribed $k$ (Examples~\ref{ex:poinc:6} and~\ref{ex:poinc:7}).

The first example of a delay of the second generation, in the author's thesis \cite[Example 5.4.8]{meng2026application}, used a rational function $R$ whose immediate basin $A_0^0$ is infinitely connected, and it was natural to suspect that delays are related to the topology of $A_0^0$. Examples~\ref{ex:poinc:6} and~\ref{ex:poinc:7} show that this is not the case: delays already occur for polynomials of every degree at least $4$, for which $A_0^0$ is simply connected. These examples are the main results of the first part of the paper.

The second part of the paper is devoted to (\ref{eq:intro:7}) and (\ref{eq:intro:8}). One main result is to show that there are no transcendental meromorphic solution of some form of FDEs. The key step in the proof relies on using Nevanlinna theory on functional equations and composition of complex functions (see \cite{goldstein1970certain,ishizaki1991conjecture,bergweiler1998meromorphic,rieppo2007class,fletcher2009meromorphic}). In our case the equation is $L[y]=P(y\circ g)$ (Theorem~\ref{thm:nev:9}): if $g$ is a polynomial of degree $\mu\geq2$, then $\tau>\sigma\mu$, where $\tau$ and $\sigma$ are the degrees of $L$ and $P$. The case $\tau=\sigma\mu$, which is exactly the case of (\ref{eq:intro:7}) and the proof requires some estimates on $A(r,y)=T(r,y)/\ln r$ along the radii $r_{i+1}=d_1r_i^\mu$ (for some constant $d_1>0$). Combined with iteration theory, this shows that the continuation cannot reach $\partial A^0_{z_0}$ on any sheet (Theorem~\ref{thm:nev:10}).

The paper is organized as follows. In Section~\ref{sec:poinc}, we recall the results in \cite[§4]{halburd2025extensionanalyticsolutionsfirstorder}, introduce stems, sheets and generations, prove the lemmas for Poincar\'e functions of higher degree, and give the examples. In Section~\ref{sec:nev}, we present Nevanlinna theory, prove the existence of solutions with a pole at $z=0$, and show that solutions of (\ref{eq:intro:7}) and (\ref{eq:intro:8}) cannot, in general, be meromorphic. Sections~\ref{sec:fde} and~\ref{sec:qde} use the language of Section~\ref{sec:poinc} and the results of Section~\ref{sec:nev} to describe the singularity structure of the solutions of (\ref{eq:intro:7}) and (\ref{eq:intro:8}). Section~\ref{sec:disc} discusses some open problems. Appendix~\ref{sec:appA} proves the existence of Laurent series solutions, and Appendix~\ref{sec:appB} contains the composition estimates used in Section~\ref{sec:nev}.

\medskip
\noindent\textbf{Notation.} Throughout, $\mathbb{N}=\{1,2,3,\dots\}$, $D(a,r)=\{z:|z-a|<r\}$ and $C(a,r)=\{z:|z-a|=r\}$. For an attracting fixed point $z_0$, $A^0_{z_0}$ denotes its immediate basin of attraction.

\section{Singularity structure of special solutions of the FE (\ref{eq:intro:3})}
\label{sec:poinc}

\subsection{Preliminaries}
\label{subsec:poinc:1}

Let $g$ be a solution of (\ref{eq:intro:1}). After the change of coordinate $w\mapsto c w$ for a suitable constant $c\neq0$, we may assume that $g$ has the expansion
\begin{equation}
    \label{eq:poinc:1}
    g(w)=c_1 w+\sum_{k\geq2}c_kw^k.
\end{equation}
Here the coefficients $c_k$ can be expressed in terms of the parameter $c_1$. A proof of the existence of such a solution can be found in the classical literature on the iteration of rational functions, see \cite[Ch.3, §4]{steinmetz2011rational} and \cite[Ch.2, §1-4]{carleson2013complex}. It is known that $g^{-1}$ maps $A_0^0$ holomorphically onto $\mathbb{C}$ \cite[Ch.3, §4 Th.2]{steinmetz2011rational}. Moreover, by the Fatou theorem \cite[Theorem 7.5.1]{beardon2000iteration}, the map $g^{-1}$ is not injective, so $R$ has a critical point $t_0\in A_0^0$. Equivalently, $g$ has a branch point. More generally, every critical point of $R$ inside $A_0^0$ gives rise to a branch point of $g$.

We now recall the result of \cite{halburd2025extensionanalyticsolutionsfirstorder}, where we described the singularity structure of the analytic continuation of $g$ for $R(t)=\lambda t+t^2$, $0<\lambda<1$. Without loss of generality, we may set $c_1=1$. The only critical point of $R$ is $t_0=-\lambda/2$, which must lie in the immediate basin of attraction $A_0^0$, and the critical value is $R(t_0)=-\lambda^2/4$. Hence, on the circle of convergence of (\ref{eq:poinc:1}), there is a branch point $w_0$ with $g(w_0)=t_0$, which results from $g(\lambda w_0)=R(t_0)$; for the position of $w_0$ see \cite[§4]{halburd2025extensionanalyticsolutionsfirstorder}. We state the following result of that work without proof:
\begin{theorem}
    \label{thm:poinc:1}
    Let $g$ be as above. Then the analytic continuation of (\ref{eq:poinc:1}) has a quadratic branch point at $w_0$. Let $g_0$ denote the function on the sheet corresponding to (\ref{eq:poinc:1}), and $g_1$ the function on the other sheet at $w_0$. Then:
    \begin{enumerate}
        \item the set of branch points of $g_0$ is $\{w_0\lambda^{1-k},k\in\mathbb{N}\}$;
        \item $g_0+g_1=-\lambda/2$, and $g_1$ has exactly the same branch points as $g_0$;
        \item each branch point of $g_0$ and $g_1$ is quadratic. At $w_0/\lambda$, there are four Riemann sheets, corresponding to the functions $g_0$, $g_1$, $g_{10}$ and $g_{11}$, which satisfy
        \begin{equation}
            g_0(\lambda w)=R(g_i(w)),\quad g_1(\lambda w)=R(g_{1i}(w)).
        \end{equation}
        \item the set of branch points of $g_{1i}$ is $\{w_0\lambda^{-k},k\in\mathbb{N}\}$;
        \item by induction, we obtain $2^{n-1}$ new Riemann sheets from the previous sheets, whose set of branch points is $\{w_0\lambda^{-k+n-1},k\in\mathbb{N}\}$. The details of how these Riemann sheets are glued are given in \cite[§4]{halburd2025extensionanalyticsolutionsfirstorder}.
    \end{enumerate}
\end{theorem}

For $R$ of higher degree, we want to show that the singularity structure can be different. To fix terminology, we introduce the following definitions for the Riemann sheets of the analytic continuation of $g$.

\medskip
\noindent\textit{Convention.} Throughout this section, a ``sheet'' means a branch $g_I$ of the analytic continuation of $g$, that is, a germ obtained by continuing $g$ from $0$ along a path and considered on a simply connected domain (a slit domain, say). Statements about ``the branch points of a sheet'' refer to the points where this branch cannot be continued analytically.

\begin{definition}[Stem, stem sheets, predecessor and successor sheets]
    \label{def:poinc:1}
    Let $R$ be the rational function (\ref{eq:intro:2}), and let $t_0,t_1,\dots,t_m$ be all its critical points. Let $g$ be the function defined by (\ref{eq:intro:1}), with $c_1=1$ in (\ref{eq:poinc:1}). Analytic continuation of $g$ gives a Riemann surface $\mathfrak{G}$ with infinitely many sheets. We label the function on each sheet by $g_I$, where $I$ is an index set, and we let $g_0$ coincide with (\ref{eq:poinc:1}). Let $\mathrm{proj}:\mathfrak{G}\to\mathbb{C}$ be the projection map from the Riemann surface to $\mathbb{C}$ given by $\mathrm{proj}((x,g_I(x)))=x$.

    Suppose that $t_0=g_j(w_0)$ and $w_0\notin \mathrm{proj}((D(0,|w_0|),g_I(D(0,|w_0|))))$. That is, the value $t_0$ is attained for the first time at $w_0$ when $g_0$ is continued analytically from $0$. Then $w_0$ is called the \emph{stem of type $t_0$}, and the sheets corresponding to $g_{j_1},\dots,g_{j_m}$, where $m$ is the multiplicity of the critical point $t_0$, are called the \emph{stem sheets of type $t_0$}. These sheets are obtained by solving $g_\ell(\lambda w)=R(g_{j_i}(w))$ for some index $\ell$. We call the sheet corresponding to $g_\ell$ the \emph{predecessor sheet} and the sheets corresponding to $g_{j_i}$ its \emph{successor sheets} of type $t_0$. We also have $g_{j_i}(w_0)=t_0$ and $g_\ell(\lambda w_0)=R(t_0)$.
\end{definition}

\begin{definition}[Branch points of a given type]
    \label{def:poinc:2}
    Assume there is a critical point $t_2$ mapping to another critical point $t_1$, that is, $R(t_2)=t_1$. The stem of type $t_1$ is called the \emph{branch point of type $t_1$}. If $w_0$ is the stem of type $t_1$, then among the successor sheets of the stem sheets of type $t_1$ there are stem sheets of type $t_2$, with stem $w_0/\lambda$. For the sheets that have a stem of type $t_2$, the branch points are called of type $t_2$ from $w_0/\lambda$ onward. For the sheets that do not have a stem of type $t_2$, the resulting branch points are called of type $t_1$. For a sheet containing both stems of type $t_1$ and of type $t_2$, only the stem $w_0$ is called a branch point of type $t_1$, and the branch points from $w_0/\lambda$ onward are of type $t_2$.

    We remark that $w_0$ can be the stem of branch points of different types.
\end{definition}

\begin{definition}[Generations]
    \label{def:poinc:3}
    We call the sheets corresponding to $g_0,\dots,g_n$ obtained by solving $g_0(\lambda w)=R(g_i(w))$ the \emph{first generation sheets}. If $g_{j_m}$ corresponds to a sheet of the $k$-th generation, all the sheets $g_{k_n}$ corresponding to solutions of $g_{j_m}(\lambda w)=R(g_{k_n}(w))$ are called the \emph{$(k+1)$-th generation sheets}.
\end{definition}

\begin{lemma}
    \label{lem:poinc:1}
        Let $t_0\in A_0^0$ be a critical point of $R$ with multiplicity $m$, i.e., $(t-t_0)^m|R(t)-R(t_0),~(t-t_0)^{m+1}\not|R(t)-R(t_0)$, and we also suppose that $t_0$ is neither a preimage nor an image of any other critical point of $R$ under the iterates $R^{\circ n}:=R\circ\cdots\circ R$ ($n$ times). Assume that $w_0$ is the stem of type $t_0$ lying on the stem sheets corresponding to $g_{j_1},\dots,g_{j_m}$. Then there are no other stems or stem sheets of type $t_0$. Moreover, if the set of branch points of type $t_0$ of these $m$ sheets is $B_0$, then the set of branch points of type $t_0$ of the successor sheets of $g_{j_i}$ is $B_0/\lambda$, and each branch point of type $t_0$ has order $m$.

\end{lemma}

\begin{proof}
    We first assume that $g_0$ is branched at $w_0$. In other words, the branch point of type $t_0$ lies on the sheet corresponding to $g_0$. Then using (\ref{eq:intro:1}), it is clear that $w_0\lambda^{-k},k\in\mathbb{N}$ are also branch points of $g_0$ of type $t_0$. If there is another branch point $w_0'\notin \{w_0\lambda^{1-k},k\in\mathbb{N}\}$ of $g_0$, then by iterating $g_0(w_0')$ with $R$, we obtain that $w_0'\lambda^n$ is another stem of type $t_0$. So $g_0(\lambda w_0)=g_0(w_0'\lambda^{n+1})=R(t_0)$ and thus, $g_0(\lambda^{k} w_0)=g_0(w_0'\lambda^{n+1+k})$. Since $w_0\lambda^k\to0$ and $g_0'(0)\neq0$, then $g_0$ is locally injective at $0$ and we obtain a contradiction. So $B_0=\{w_0\lambda^{1-k},k\in\mathbb{N}\}$. By counting multiplicities of the critical point $t_0$, solving (\ref{eq:intro:1}) yields exactly $m$-th sheets with the stem $w_0$ of type $t_0$ whose set of branch points of type $t_0$ is also $B_0$. If $h$ is a succesor sheet of $g_i$, then it is clear that $h$ is branched at $\{w_0\lambda^{-k},k\in\mathbb{N}\} $. If $w_0$ is a branch point of type $t_0$ of $h$, then either $\lambda w_0$ is also branch point of type $t_0$ of $g_i$, the predecessor of $h$, which contradicts the result of $B_0$. Otherwise $w_0$ is the stem of type $t_0$ of $h$, leading to $g_i(\lambda w_0)=R(t_0)$. However, $g_0(\lambda w_0)=R(t_0)$ and $R(t_0)$ is not a critical point of $R$. Thus, it is impossible that $g_i(\lambda w_0)=g_0(\lambda w_0)$ since both arised from solving $g_0(\lambda w)=R(g_j(w))$. Hence, $B_0/\lambda$ is the set of branch points of type $t_0$ of $h$. In the case where the stem lies on the $k$-th generation sheet corresponding to $g_j\neq g_0$, we iterate the point $w_0$ to $\lambda^kw_0$. From there, the argument for $g_0$ and $g_i$ can be used to derive the conclusion in the same way.

    To show that the branch points of type $t_0$ have order $m$, it is enough to consider the following factorisation:
    \[
    g(\lambda w)-R(t_0)=(g(w)-t_0)^mR_1(g(w)),
    \]
    where $R_1$ is the remaining term from factoring out all factor $t-t_0$. Then the Puiseux series of $g(w)$  is of the form
    \[
    g(w)=t_0+C_1\sqrt[m]{g(\lambda w)-R(t_0)}+\cdots=t_0+C_1'\sqrt[m]{g'(\lambda w_0)(w-w_0)+\cdots}+\cdots.
    \]
    If $g_0$ is branched on $w_0$, it is enough to show that $g_0'(\lambda w_0)\neq0$. Assume the contrary. Differentiating (\ref{eq:intro:1}) and substituting $w=\lambda^kw_0$ gives $g_0'(\lambda^kw_0)=0$ for all $k\geq1$. Since $\lambda^kw_0\to0$, the zeros of $g_0'$ accumulate at $0$, so $g_0'\equiv0$ by the identity theorem, contradicting $g_0'(0)=1$.
    If the stem $w_0$ does not lie on the first generation, say $g_{k_n}$ and assume that $g_{k_{n-1}}$ is its predecesor, then we obtain
    \[
    g_{k_n}(w)=t_0+C_1\sqrt[m]{g_{k_{n-1}}(\lambda w)-R(t_0)}+\cdots=t_0+C_1'\sqrt[m]{g_{k_{n-1}}'(\lambda w_0)(w-w_0)+\cdots}+\cdots.
    \]
    If $g_{k_{n-1}}'(\lambda w_0)=0$, then (\ref{eq:intro:1}) yields $g_{k_{n-2}}'(\lambda^2w_0)=0$, where $g_{k_{n-2}}$ is the predecessor of $g_{k_{n-1}}$. Thus, for $k$ large enough, we obtain $g_0'(\lambda^kw_0)=0$, and a contradiction follows as above.
\end{proof}

\begin{lemma}
    \label{lem:poinc:2}
    Suppose that there is a critical point $t_2\in A_0^0$ with multiplicity $m_2$ mapping to another critical point $t_1\in A_0^0$ with multiplicity $m_1$, i.e., $R(t_2)=t_1$. Suppose further that no other critical point of $R$ lies in $R^{\pm\circ n}(t_1)$ for any $n$.
    We assume that $t_2$ is not an image of any other critical point of $R$. Let $w_1$ be the stem of type $t_1$ lying on the sheets corresponding to $g_{j_1},\dots,g_{j_{m_1}}$, then there are $m_1m_2$ stem sheets of type $t_2$ which are the successors of the sheets corresponding to $g_{j_1},\dots,g_{j_{m_1}}$ with the stem $w_1/\lambda$. We denote the functions corresponding to those stem sheets by $h_{1},\dots,h_{m_1m_2}$. Let $B_1$ be the set of branch points of type $t_1$ of $g_{j_i}$, and $B_2,B_2'$ be the sets of branch points of type $t_2$ of $g_{j_i}$ and $h_k$ respectively. Then the set of branch points of type $t_i$ of the successor sheets of $g_{j_i}$, other than the $h_k$, is $B_i/\lambda$, and each such branch point has order $m_1$. For the $h_k$, we obtain $B_2'=B_2/\lambda$, and each branch point has order $m_1m_2$.

    If $R^{\circ n}(t_2)=t_1$, the same arguments apply, starting from $w_1/\lambda^n$ instead of $w_1/\lambda$.

    If $t_2$ is the image of another critical point, one accumulates the orders starting from the stem $w_1/\lambda^2$; we do not comment on this further, as it follows by induction.
\end{lemma}

\begin{proof}
    The statements about the sets of branch points of types $t_1$ and $t_2$ are proved in the same way as in Lemma \ref{lem:poinc:1}, and we do not repeat the details.
    For the order of the branch points of type $t_2$, note that
    \begin{equation}
    \begin{aligned}
    g(\lambda^2w)-R(t_1)&=g(\lambda^2w)-R(R(t_2))=(g(\lambda w)-t_1)^{m_1}R_1(g(\lambda w))\\
    &=(g(\lambda w)-R(t_2))^{m_1}R_1(g(\lambda w)),\\
    g(\lambda w)-R(t_2)&=g(\lambda w)-t_1=(g(w)-t_2)^{m_2}R_2(g(w)).
    \end{aligned}
    \end{equation}
    Solving for $g(w)$, we see that the branch point has order $m_1m_2$ in terms of $g(\lambda^2w)-R(t_1)$, and we argue as in Lemma \ref{lem:poinc:1} to show that $g'(\lambda w_0)\neq0$.
\end{proof}

\begin{lemma}
    \label{lem:poinc:3}
    Suppose that the hypotheses of Lemmas \ref{lem:poinc:1} and \ref{lem:poinc:2} hold, and let $w_0$ be the stem of type $t_0$.
    \begin{enumerate}
        \item If $g_0$ is branched at $w_0$, then the set of branch points of type $t_0$ on $g_0$ is $B_0=\{w_0\lambda^{1-k},k\in\mathbb{N}\}$, and on the sheets of the $n$-th generation it is $B_0\lambda^{1-n}$.
        \item If $g_0$ is not branched at $w_0$, then on each sheet containing the stem $w_0$, the only branch point of type $t_0$ is $w_0$ itself.
    \end{enumerate}
\end{lemma}

\begin{proof}
    Case (1) follows from the proof of Lemma \ref{lem:poinc:1}. For case (2), suppose first that the stem $w_0$ lies on another first generation sheet, say $g_1$. We claim that $w_0\lambda^{-k}$, $k\in\mathbb{N}$, is not a branch point of $g_1$. Otherwise, $R(t_0)$ would be a value of $g_0$ at a point $w_0\lambda^{1-k}$ with $k\geq1$. Since also $g_0(\lambda w_0)=R(t_0)$, the local injectivity argument in the proof of Lemma \ref{lem:poinc:1} gives a contradiction.

    Suppose next that $w_0$ lies on a sheet of a later generation, say $h_1$, with predecessor $h_0$, and that $w_0\lambda^{-k}$ is a branch point of $h_1$ for some $k$. Then again $R(t_0)$ is a value of $h_0$ at $w_0\lambda^{1-k}$, and also $h_0(\lambda w_0)=R(t_0)$. Applying (\ref{eq:intro:1}) repeatedly, we arrive at $g_0(\lambda^nw_0)=g_1(w_0\lambda^{n+1-k})$ for some $n$. Applying (\ref{eq:intro:1}) once more shows that, if $k\neq1$, the sequence $g_0(\lambda^{n+mk}w_0)$, $m\geq0$, is constant; since $\lambda^{n+mk}w_0\to0$, this forces $g_0$ to be constant, a contradiction. Hence $k=1$. But then $g_0(\lambda^nw_0)=g_1(\lambda^nw_0)$, so $\lambda^nw_0$ would itself be a branch point, which contradicts the assumption that no other critical point of $R$ lies in the orbit $R^{\pm\circ n}(t_0)$. Therefore $w_0$ is the only branch point of type $t_0$ on these sheets.
\end{proof}

\begin{lemma}
    \label{lem:poinc:4}
Let $h_i$ and $h_{i+1}$ be the functions on two different sheets with the same predecessor $g_{i-1}$. Then $h_i$ and $h_{i+1}$ cannot share more than one common value. The common value is attained at a branch point of the same type lying at the stem $w_0$.
\end{lemma}

\begin{proof}
    Assume that $h_i(v_1)=h_{i+1}(v_0)$ such that $v_1\neq v_0$. Then, applying (\ref{eq:intro:1}) repeatedly to pass down to the sheet corresponding to $g_0$, we obtain $g_0(\lambda^k v_1)=g_0(\lambda^k v_0)$ and the local injectivity at $0$ proves that $v_1=v_0$. Since $v_1=v_0$, then the image of $h_i(v_1)=h_{i+1}(v_0)$ is a critical point and so $v_1=v_0$ is also a stem. Suppose that there is another point $v_0'\neq v_0$ such that $h_i(v_0')=h_{i+1}(v_0')$. Near the stem $v_0$, we have, for $h\in\{h_i,h_{i+1}\}$, the factorisation
    \begin{equation}
        g_{i-1}(\lambda w)-R(h_i(v_0))=(h(w)-h_i(v_0))^mR_1(h(w)),
    \end{equation}
    where $R_1$ does not vanish at $h_i(v_0)$ or at $h_i(v_0')$. Taking $m$-th roots, we obtain
    \begin{equation}
        \begin{aligned}
        &[h_i(w)-h_i(v_0)]\sqrt[m]{R_1(h_i(w))}=e^{\frac{2\pi ik_1}{m}}\Phi(w),\\
        &[h_{i+1}(w)-h_i(v_0)]\sqrt[m]{R_1(h_{i+1}(w))}=e^{\frac{2\pi ik_2}{m}}\Phi(w),
        \end{aligned}
    \end{equation}
    where $\Phi(w)=\sqrt[m]{g_{i-1}(\lambda w)-R(h_i(v_0))}$ and $k_1\neq k_2$, since $h_i\neq h_{i+1}$. Eliminating $\Phi$, we get
    \begin{equation}
        e^{\frac{2\pi ik_2}{m}}[h_i(w)-h_i(v_0)]\sqrt[m]{R_1(h_i(w))}-e^{\frac{2\pi ik_1}{m}}[h_{i+1}(w)-h_i(v_0)]\sqrt[m]{R_1(h_{i+1}(w))}=0.
    \end{equation}
    We continue this identity analytically along a path from $v_0$ to $v_0'$ on which $h_i$ and $h_{i+1}$ are analytic, choosing the branches of the $m$-th roots continuously along the path. Setting $w=v_0'$ and using $h_i(v_0')=h_{i+1}(v_0')$, we obtain
    \begin{equation}
        (e^{\frac{2\pi ik_2}{m}}-e^{\frac{2\pi ik_1}{m}})[h_i(v_0')-h_i(v_0)]\sqrt[m]{R_1(h_i(v_0'))}=0,
    \end{equation}
    which is possible only when $h_i(v_0')=h_i(v_0)$.
\end{proof}

\subsection{Examples of singularity structure}
\label{subsec:poinc:2}

When $R$ is a quadratic polynomial, there is only one critical point, and the stem of the branch point must lie on the sheet corresponding to $g_0$ \cite[§4]{halburd2025extensionanalyticsolutionsfirstorder}. In this case, the singularity structure of all the Riemann sheets of the Poincar\'e function can be described as in Lemmas \ref{lem:poinc:1} and \ref{lem:poinc:3}(1). In the examples below, we illustrate different singularity structures of the Poincar\'e function when $R$ has higher degree, in particular the situation of Lemma \ref{lem:poinc:3}(2), where the stem can lie on a sheet of any generation.

We start with the simplest example, in which all stems lie on the sheet corresponding to $g_0$. Here we normalise (\ref{eq:poinc:1}) so that it matches the Poincar\'e function (\ref{eq:intro:1}), that is, $c_1=1$. In all the examples below, $R$ has real coefficients, so within the radius of convergence of (\ref{eq:poinc:1}) the standard results of real analysis can be applied. 
The facts that the critical points lie in $A_0^0$, and the order in which the critical values are attained, can be checked numerically or plotting the Julia set of $R$.

\begin{example}
\label{ex:poinc:1}
Consider the equation
\begin{equation}
\label{eq:poinc:2}
    g(\lambda^2w)=R_1(g(w))=\lambda^2g(w)-\frac{\lambda^4}{27}g(w)^3,\quad\lambda\in(0,1).
\end{equation}
The function $R_1$ admits two critical points $\pm3/\lambda$, both are in $A_0^0$. We set $t_0=3/\lambda$ and $t_0'=-3/\lambda$. Moreover, $R_1(\pm3/\lambda)=\pm2\lambda$, and
\begin{subequations}
\label{eq:poinc:3}
\begin{align}
\label{eq:poinc:4}
    &g(\lambda^2w)-2\lambda=-\frac{\lambda^4}{27}\left(g(w)-\frac{3}{\lambda}\right)^2\left(g(w)+\frac{6}{\lambda}\right),\\
\label{eq:poinc:5}
    &g(\lambda^2w)+2\lambda=-\frac{\lambda^4}{27}\left(g(w)+\frac{3}{\lambda}\right)^2\left(g(w)-\frac{6}{\lambda}\right).
\end{align}
\end{subequations}
One checks that the power series of $g$ in (\ref{eq:poinc:1}) is odd, that is, $c_{2k}=0$. Since the continuation of the series must have a branch point $w_0$, the point $-w_0$ is also a branch point. Since the stems of both branch points lie on the sheet corresponding to $g_0$, the singularity structure of each type of branch point corresponds to case (1) of Lemma \ref{lem:poinc:3}, and all other Riemann sheets have their singularities distributed as described in Lemma \ref{lem:poinc:1}. The branch point structure is illustrated in Figure \ref{fig:poinc:1}(A). The black line corresponds to $g_0$, the orange line corresponds to $g_1$, which shares the stem of type $t_0'$ (red dots), and the purple line corresponds to $g_2$, which shares the stem of type $t_0$ (blue dots). We see from the picture that $g_1(w_0)=-6/\lambda$ on the orange line and, similarly, $g_2(-w_0)=6/\lambda$ on the purple line; neither $g_1$ at $w_0$ nor $g_2$ at $-w_0$ is branched. The red dots on each line represent the branch points $-w_0/\lambda^k$ of type $t_0'$, and similarly for the blue dots. In Figure \ref{fig:poinc:1}(B), there are six additional sheets, obtained as the successors of $g_1$ and $g_2$. We see how the branch points of each type are distributed on these sheets, following their predecessors. For example, the yellow line corresponds to a successor sheet of $g_1$, and its branch points are shifted by a factor $1/\lambda$ from those of $g_1$, so there is one fewer branch point of each type. We do not show how these sheets are glued; this requires further analysis as in \cite{halburd2025extensionanalyticsolutionsfirstorder} and is not discussed in this paper.

\figPoincareOne

\end{example}

For the next few examples we omit the details of the analysis, as it is the same as in Example \ref{ex:poinc:1}. The following example illustrates Lemma \ref{lem:poinc:2}.

\begin{example}
\label{ex:poinc:2}
Consider the equation
\begin{equation}
\label{eq:poinc:6}
g\left(\lambda w\right)=R_2(g(w))=\frac{-3}{20}g(w)+\frac{-12}{20}g^2(w)+g^3(w),~\lambda=\frac{-3}{20}.
\end{equation}
Here, we have the following relations:
\begin{subequations}
\label{eq:poinc:7}
\begin{align}
\label{eq:poinc:8}
&g\left(\lambda w\right)-\frac{1}{125}=\left(g(w)-\frac{-1}{10}\right)^2\left(g(w)-\frac{4}{5}\right),\\
\label{eq:poinc:9}
&g\left(\lambda w\right)-\frac{-1}{10}=\left(g(w)-\frac{1}{2}\right)^2\left(g(w)-\frac{-2}{5}\right).
\end{align}
\end{subequations}
The critical points $-1/10$ and $1/2$ both lie in $A_0^0$ and $R_2(1/2)=-1/10$. We set $t_1=-1/10$ and $t_2=1/2$. Since $R_2(-1/10)=1/125$, the function $g_0$ will radially hit the value $1/125$ first. The value $1/2$ cannot be attained first, since this would require $g_0$ to attain $R_2(1/2)=-1/10$ at a point of smaller modulus.  We then obtain the branch point $w_0$ of type $t_1$ and order $2$, which is shared with the sheet corresponding to $g_1$, represented by the black and orange lines respectively in Figure \ref{fig:poinc:2}(A). By Lemma \ref{lem:poinc:4}, we see that $g_0$ and $g_1$ share only one common value. Thus, lying on the stem $w_2$ of type $t_2$, $g_1(w_2)=-2/5$. This leaves another successor sheet obtained from  $R_2^{-1}g_0$ to have a branch point at $w_2$ which is a stem of type $t_2$. The branch point $w_2$ of $g_0$ and $g_2$ is of fourth order and the other two sheets are obtained as successors of $g_1$, represented by the light green and yellow lines in Figure \ref{fig:poinc:2}(B). The branch points of the second generation sheets follow from Lemma \ref{lem:poinc:2}.

\figPoincareTwo

\end{example}

The situation of Lemma \ref{lem:poinc:3}(2) does not occur in Examples \ref{ex:poinc:1} and \ref{ex:poinc:2}; it occurs in the next example.

\begin{example}
\label{ex:poinc:3}
Consider the equation
\begin{equation}
\label{eq:poinc:10}
g\left(\lambda w\right)=R_3(g(w))=\frac{2}{3}g(w)-\frac{3}{2}g(w)^2+g(w)^3,~\lambda=\frac{2}{3}.
\end{equation}
Here $R_3$ has two critical points, $t_0=1/3$ and $t_0'=2/3$, both of which are quadratic and lie in $A_0^0$.  We also have
\begin{subequations}
\label{eq:poinc:11}
\begin{align}
\label{eq:poinc:12}
&g(\lambda w)-\frac{2}{27}=\left(g(w)-\frac{1}{6}\right)\left(g(w)-\frac{2}{3}\right)^2,\\
\label{eq:poinc:13}
&g(\lambda w)-\frac{5}{54}=\left(g(w)-\frac{1}{3}\right)^2\left(g(w)-\frac{5}{6}\right).
\end{align}
\end{subequations}
Arranging the values $2/27,5/54,1/6,2/3,1/3,5/6$ on the real line, we obtain from the Intermediate Value Theorem (IVT) that $g_0$ radially hits the value $2/27$ first, at $\lambda w_0$, and then $5/54$, at $\lambda w_0'$; see the black line in Figure \ref{fig:poinc:3}(A). On $g_0$, the point $w_0$ cannot be a stem of type $t_0'$: otherwise $g_0(w_0)=2/3>1/3$, and by the IVT the value $1/3$ would have to appear on $g_0$ before $2/3$, which is impossible since $R_3(2/3)=2/27<R_3(1/3)=5/54$. Hence $g_0(w_0)=1/6$, which is a regular value, and the first branch point on $g_0$ is $w_0'$, with $g_0(w_0')=1/3$, of type $t_0$. The value $2/3$, a double root of $R_3(t)=2/27$, is taken at $w_0$ by the two other first generation sheets $g_1$ and $g_2$, obtained by solving $g_0(\lambda w)=R_3(g_i(w))$ and represented by the orange and purple lines in Figure \ref{fig:poinc:3}(A). Thus $w_0$ is the stem of type $t_0'$, and by Lemma \ref{lem:poinc:3}(2) it is the only branch point of type $t_0'$ on $g_1$ and $g_2$. In the second generation, the successor sheets of $g_1$ and $g_2$ have the unique branch point $w_0/\lambda$ of type $t_0'$, as seen on the additional six lines (except the pink line) in Figure \ref{fig:poinc:3}(B); the pink line represents a third generation sheet whose predecessor is represented by either the light green or the lime green line.  For type $t_0$, the distribution of the branch points is the same as in the previous Examples and Lemma \ref{lem:poinc:1}.

This example shows that the type $t_0'$ has only one element on each Riemann sheet (except the sheet corresponding to $g_0$).

\figPoincareThree

\end{example}

The next example shows that the stem can be at $0$; this occurs when a critical point is a root of $R$.

\begin{example}
\label{ex:poinc:4}
Consider the equation
\begin{equation}
\label{eq:poinc:14}
g(\lambda w)=R_4(g(w))=g(w)\left(g(w)-\frac{1}{3}\right)^2,~\lambda=\frac{1}{9}.
\end{equation}
Here, $R_4$ has two critical points $t_0=1/9$ and $t_0'=1/3$. We also have the relation
\begin{equation*}
\label{eq:poinc:15}
g(\lambda w)-\frac{4}{729}=\left(g(w)-\frac{1}{9}\right)^2\left(g(w)-\frac{4}{9}\right).
\end{equation*}
Here the distribution of the branch points is immediate, since the critical point $1/3$ must be the image of the branch point $0$ lying on the first generation sheets other than $g_0$; see the orange and purple lines in Figure \ref{fig:poinc:4}. The branch point of type $t_0$ has to appear on $g_0$ and is glued with $g_1$ at the stem $w_0$.

\figPoincareFour

\end{example}

It is also possible that $w_0$ is the stem of two different branch points lying on different sheets, as in the following example.

\begin{example}
\label{ex:poinc:5}
Consider the equation
\begin{equation}
\label{eq:poinc:16}
g(\lambda w)=R_5(g(w))=\frac{-5}{18}g(w)+\frac{37}{36}g^2(w)-\frac{5}{3}g^3(w)+g^4(w),~\lambda=\frac{-5}{18}.
\end{equation}
We also have
\begin{align*}
&g(\lambda w)-\left(\frac{-575}{20736}\right)=\left(g(w)-\frac{5}{12}\right)^2\left(g^2(w)-\frac{5}{6}g(w)+\frac{23}{144}\right),   \\
&g(\lambda w)-\left(\frac{-1}{36}\right)=\left(g(w)-\frac{1}{3}\right)^2\left(g(w)-\frac{1}{2}\right)^2.
\end{align*}
Here $R_5$ has three critical points, $t_0=5/12$, $t_0'=1/3$ and $t_0''=1/2$, all lying in $A_0^0$.
Using the IVT as in Example \ref{ex:poinc:3}, we find that $g_0$ hits the value $1/6$ and that the stem $w_0$ is a branch point of type $t_0$ of $g_1$ and $g_2$ (orange and purple lines in Figure \ref{fig:poinc:5}). For the sheet with the value $g_0$, $g_0$ is branched at $w_0'=w_0''$ such that $g_0(w_0')=1/3$ by IVT. Using Lemma \ref{lem:poinc:4}, we can see that the sheets corresponding to $g_2$ and $g_3$ (brown line in Figure \ref{fig:poinc:5}) are branched of type $t_0''$ with the stem $w_0'=w_0''$. Thus, $g_1$ is branched of type $t_0'$ at $w_0'$; the matching here requires Lemma \ref{lem:poinc:4}.

\figPoincareFive

\end{example}

In \cite[Example 5.4.8, p.171]{meng2026application}, we constructed an example with a rational function $R$ such that the stem of a branch point appears on a second generation sheet. The construction relies on the following ideas. Let $t_0$ and $t_0'$ be two critical points of $R$. If $g_0$ is radially branched first of type $t_0$ at $w_0$ such that
\begin{equation}
    \label{eq:poinc:17}
    R^{\circ2}(t_0')<t_0<R(t_0'),
\end{equation}
then we may use the Puiseux series of $g_0$ at $w_0/\lambda$, or the preimage values of $R^{-1}(R^{\circ2}(t_0'))$, to show that the value $R(t_0')$ cannot be a value of $g_0$. Thus it must be a value of another first generation sheet, say $g_1$. So one of the successors of $g_1$ takes the value $t_0$, where the function is branched at the stem $w_0/\lambda^2$ of type $t_0'$.

In the same spirit, we construct a similar situation for a polynomial $P$. In this case, the form of (\ref{eq:intro:1}) is
\begin{equation}
    \label{eq:poinc:18}
    g(\lambda w)=0.51\int_0^{g(w)}(1-u)^2(1-\frac{u}{4})du.
\end{equation}
Using this ansatz, we prove a generic case in which the stem first occurs on a sheet of the $k$-th generation, for some polynomial of even degree. The case (\ref{eq:poinc:18}) is covered by Example \ref{ex:poinc:6}, where the constant $0.51$ is also explained.

\begin{example}
\label{ex:poinc:6}
\upshape
Let $n\geq2$ be an even natural number, and let $k\geq2$ be a natural number. (For $k=1$, condition (2) below reads $R_{\lambda,r}(r)<1<r$, which holds for every sufficiently small $\lambda>0$.)
Consider the equation
    \begin{equation}
    \label{eq:poinc:19}
    g(\lambda w)=R_{\lambda,r}(g(w))=\lambda\int_0^{g(w)}(1-u)^{n}\left(1-\frac{u}{r}\right)du,\quad r>1.
    \end{equation}
    The polynomial $R_{\lambda,r}$ has two critical points, $1$ and $r$, and $R_{\lambda,r}'(t)\geq0$ for all $t\in[0,r]$.
    We show that, for every sufficiently large $r$ (depending on $k$), there exists $0<\lambda<1$ such that
    \begin{enumerate}
        \item $R_{\lambda,r}([0,r])\subset[0,r]$, that is, the interval $[0,r]$ is forward invariant, and so it lies in $A_0^0$;
        \item $R_{\lambda,r}^{\circ k}(r)<1<R_{\lambda,r}^{\circ (k-1)}(r)$;
    \end{enumerate}
    and, moreover, that the value $R_{\lambda,r}^{\circ (k-1)}(r)$ is not attained by $g_0$. We write
    \[
    S_{r}(t)=\int_0^t(1-u)^{n}\left(1-\frac{u}{r}\right)du,
    \]
    so that $R_{\lambda,r}=\lambda S_r$. The function $S_r$ is non-decreasing on $[0,r]$, and
    \begin{equation}
        \label{eq:poinc:20}
        S_r(r)=\frac{1}{n+1}-\frac{1}{r(n+1)(n+2)}+\frac{(1-r)^{n+2}}{r(n+1)(n+2)},
    \end{equation}
    so $S_r(r)$ is of order $r^{n+1}$ as $r\to\infty$.

    We first show that there exists $\lambda$ such that condition (2) holds. Fix a large $r$, and let $\lambda_0=r/S_r(r)\in(0,1)$, so that $r$ is a fixed point of $R_{\lambda_0,r}$:
    \begin{equation}
        \label{eq:poinc:21}
        R_{\lambda_0,r}(r)=r.
    \end{equation}
    Then $R_{\lambda_0,r}^{\circ(k-1)}(r)=r>1$, while $R_{0,r}^{\circ(k-1)}(r)=0$ since $R_{0,r}\equiv0$. By the intermediate value theorem (IVT), applied to $R_{\lambda,r}^{\circ(k-1)}(r)$ as a function of $\lambda$, there exists $\lambda\in(0,\lambda_0)$ with $R_{\lambda,r}^{\circ(k-1)}(r)=1$; let $\lambda_1$ be the smallest such value. Then
    \[
    R_{\lambda_1,r}^{\circ k}(r)=R_{\lambda_1,r}(1)=\lambda_1S_r(1)<\frac{\lambda_1}{n+1}<1.
    \]
    Let
    \begin{equation}
        \label{eq:poinc:23}
        \lambda_3=\inf_{t\in(0,r]}\frac{t}{S_r(t)}.
    \end{equation}
    Since $t/S_r(t)$ is continuous and positive on $(0,r]$ and tends to $1$ as $t\to0^+$, we have $0<\lambda_3\leq1$. We will show below that
    \begin{equation}
        \label{eq:poinc:24}
        R_{\lambda_3,r}^{\circ(k-1)}(r)>1.
    \end{equation}
    Assume that (\ref{eq:poinc:24}) holds. Since $\lambda_1$ is the smallest root of $R_{\lambda,r}^{\circ(k-1)}(r)-1$, and $R_{\lambda,r}^{\circ(k-1)}(r)$ is continuous in $\lambda$ and vanishes at $\lambda=0$, then we obtain $\lambda_1<\lambda_3$. Thus, we may choose $\lambda_2\in(\lambda_1,\lambda_1+\epsilon)\subset(\lambda_1,\lambda_3)$ with small enough $\epsilon>0$, such that
    \begin{equation}
        \label{eq:poinc:22}
        R_{\lambda_2,r}^{\circ k}(r)<1=R_{\lambda_1,r}^{\circ (k-1)}(r)<R_{\lambda_2,r}^{\circ (k-1)}(r).
    \end{equation}
    Thus, we prove condition (2); here, we can make $R_{\lambda_2,r}^{\circ (k-1)}(r)$ as close to $1$ as we wish.
    
    Since $\lambda_2<\lambda_3$, we find that $R_{\lambda_2,r}(t)=\lambda_2S_r(t)<\lambda_3S_r(t)\leq t$ for $t\in(0,r]$, which gives (1). 
    So we only need to show that (\ref{eq:poinc:24}) holds for large $r$. Substituting $u=rv$, we have
    \[
    S_r(x)=\int_0^{x/r}(1-rv)^n(1-v)\,r\,dv.
    \]
    Since $n$ is fixed and $x/r\leq1$, for all large $r$ we have $(1-rv)^n\asymp r^nv^n$, and thus $S_r(x)\asymp r^{n+1}\int_0^{x/r}v^n(1-v)\,dv$. To justify this, note that for small $x$ the function $S_r$ is differentiable at $x=0$ with $S_r'(0)=1\neq0$, so $S_r(x)\sim x$ as $x\to0^+$; and for $x/r\in[0,1]$ away from $0$, the integral $\int_0^{x/r}v^n(1-v)\,dv$ is non-zero. The asymptotics of $S_r(x)$ follow from these two facts. Hence
    \[
    \lambda_3=\inf_{x\in(0,r]}\frac{x}{S_r(x)}\asymp\inf_{x\in(0,r]}\frac{1}{r^n}\cdot\frac{x/r}{\int_0^{x/r}v^n(1-v)\,dv}.
    \]
    The function $s\mapsto s/\int_0^{s}v^n(1-v)\,dv$ is continuous and positive on $(0,1]$ and tends to $1$ as $s\to0^+$, so it has a positive minimum. Hence
    \[
    \lambda_3\in\left(\frac{a_1}{r^n},\frac{a_2}{r^n}\right)
    \]
    for some constants $0<a_1<a_2$ depending only on $n$. Since also $S_r(r)/r^{n+1}\asymp\int_0^1v^n(1-v)\,dv$, we obtain
    \[
    \frac{R_{\lambda_3,r}(r)}{r}=\frac{\lambda_3S_r(r)}{r}\geq u_1
    \]
    for some constant $u_1>0$ independent of $r$. So $R_{\lambda_3,r}(r)\geq u_1r>1$ for all large $r$. Moreover, $R_{\lambda_3,r}(r)=\lambda_3S_r(r)\leq r$, so $u_1<1$, and we may use the monotonicity of $R_{\lambda_3,r}$ and $S_r$ on $[0,r]$. For the second iterate,
    \[
    \frac{R_{\lambda_3,r}^{\circ2}(r)}{r}\geq\frac{\lambda_3S_r(u_1r)}{r}\geq u_2>0,
    \]
    since $S_r(u_1r)/r^{n+1}\asymp \int_0^{u_1}v^n(1-v)\,dv$ has a positive lower bound for all large $r$. So $R_{\lambda_3,r}^{\circ2}(r)>u_2r>1$ for all large $r$, and $u_2<1$ because $R_{\lambda_3,r}^{\circ2}(r)\leq R_{\lambda_3,r}(r)\leq r$. By induction, for each $j$ there is $0<u_j<1$ such that $R_{\lambda_3,r}^{\circ j}(r)>u_jr>1$ for all large $r$. Taking $j=k-1$ gives (\ref{eq:poinc:24}).

    Now we need to show that the value $R_{\lambda_2,r}^{\circ (k-1)}(r)$ is not attained by $g_0$. Let $w_0$ be the stem of type $1$, that is, $g_0(w_0)=1$, so that $g_0(\lambda_2w_0)=R_{\lambda_2,r}(1)$, and $g_0$ is real and increasing on $[0,w_0)$. Since $R_{\lambda_2,r}^{\circ(k-1)}(r)>1$ and $R_{\lambda_2,r}$ is increasing on $[0,r]$, we have $R_{\lambda_2,r}(1)<R_{\lambda_2,r}^{\circ k}(r)<1$. Hence there is $\tilde{w}_0>w_0$, with $\lambda_2\tilde{w}_0$ inside the disc of convergence of (\ref{eq:poinc:1}), such that $g_0(\lambda_2\tilde{w}_0)=R_{\lambda_2,r}^{\circ k}(r)$. Moreover, $\tilde{w}_0$ is close to $w_0$, since $R_{\lambda_2,r}^{\circ(k-1)}(r)$ is close to $1$ by the choice of $\lambda_2$.
    We also have
    \begin{equation}
        \label{eq:poinc:25}
    R_{\lambda_2,r}(1+x)-R_{\lambda_2,r}(1)=D_1x^{n+1}(D_2-x),
    \end{equation}
    where
    \[
    D_1=\frac{\lambda_2}{r(n+2)},\quad D_2=\frac{(n+2)(r-1)}{n+1}.
    \]
    Solving (\ref{eq:intro:1}) near $w_0$, we find that
    \begin{equation}
        \label{eq:poinc:26}
        g_0(w)=1+e^{2\pi ij/(n+1)}\left(\frac{g_0(\lambda_2 w)-R_{\lambda_2,r}(1)}{D_1D_2}\right)^{1/(n+1)}+\cdots
    \end{equation}
    for some $j\in\{0,\dots,n\}$, where $(\cdot)^{1/(n+1)}$ denotes the principal branch. We want to show that $j\neq0$. Then, for real $w>w_0$ close to $w_0$, the quantity in the brackets is positive, so $g_0(w)$ is not real, and the real values slightly larger than $1$, in particular $R_{\lambda_2,r}^{\circ (k-1)}(r)$, are taken on the branch corresponding to $j=0$.

    To see that $j\neq0$, first note that $g_0(w)<1$ and $g_0(\lambda_2w)<R_{\lambda_2,r}(1)$ for $0<w<w_0$. We set $x=g_0(w)-1<0$ and $G=R_{\lambda_2,r}(1)-g_0(\lambda_2 w)>0$, so that
    \begin{equation}
        \label{eq:poinc:27}
    R_{\lambda_2,r}(1+x)-R_{\lambda_2,r}(1)=g_0(\lambda_2 w)-R_{\lambda_2,r}(1)=-G,
    \end{equation}
    that is, $D_1x^{n+1}(D_2-x)=-G$.
    We change to the positive variable $X=-x$. Since $n+1$ is odd,
    \begin{equation}
        \label{eq:poinc:28}
    G=D_1X^{n+1}(D_2+X)\implies X=\left(\frac{G}{D_1D_2}\right)^{1/(n+1)}+O\big(G^{2/(n+1)}\big),
    \end{equation}
    and thus
    \begin{equation}
        \label{eq:poinc:29}
    x=-X=-\left(\frac{G}{D_1D_2}\right)^{1/(n+1)}+O\big(G^{2/(n+1)}\big).
    \end{equation}
    Substituting back $x$ and $G$, we obtain
    \begin{equation}
        \label{eq:poinc:30}
    g_0(w)=1-\left(\frac{R_{\lambda_2,r}(1)-g_0(\lambda_2 w)}{D_1D_2}\right)^{1/(n+1)}+\cdots,
    \end{equation}
    with the positive real root. For $0<w<w_0$, the quantity in the brackets in (\ref{eq:poinc:26}) is negative, so its principal root is $e^{i\pi/(n+1)}$ times the positive root in (\ref{eq:poinc:30}). Comparing (\ref{eq:poinc:26}) and (\ref{eq:poinc:30}) gives $e^{i\pi(2j+1)/(n+1)}=-1$, that is, $j=n/2$.

    Finally, to show that there is no other point $w$ with $g_0(w)=R_{\lambda_2,r}^{\circ (k-1)}(r)$, it suffices to show that $\lambda_2\tilde{w}_0$ is the only point at which $g_0$ takes the value $R_{\lambda_2,r}^{\circ k}(r)$. Indeed, if $g_0(w)=R_{\lambda_2,r}^{\circ (k-1)}(r)$, then $g_0(\lambda_2 w)=R_{\lambda_2,r}^{\circ k}(r)$ by (\ref{eq:intro:1}), so $w=\tilde{w}_0$; but $g_0(\tilde w_0)$ is not real, since $j\neq0$ and $\tilde w_0$ is close to $w_0$. The uniqueness of $\lambda_2\tilde{w}_0$ follows from the local injectivity argument in the proof of Lemma \ref{lem:poinc:1}, and we omit the details.

    Example (\ref{eq:poinc:18}) correspond to$n=2$ and $r=4$, and a direct computation gives $S_4(4)=\int_0^4(1-u)^2(1-u/4)\,du=2$. So $\lambda=0.5$ gives $R_{0.5,4}(4)=1$, and the slightly larger value $\lambda=0.51$ gives $R_{0.51,4}(4)=1.02>1$. Moreover, $R_{0.51,4}^{\circ2}(4)=0.51\,S_4(1.02)\approx0.159<1$, so this is the case $k=2$. Since $(1-u)^2(1-u/4)\leq1$ on $[0,4]$, we have $S_4(t)\leq t$, hence $\inf_{t\in(0,4]}t/S_4(t)=1>0.51$, and $[0,4]$ is forward invariant. These values can also be checked numerically.
\end{example}

We now extend the construction to polynomials of odd degree, by adding a further factor whose zero, a critical point, lies in the basin of infinity.

\begin{example}
    \label{ex:poinc:7}
    Consider
    \begin{equation}
    \label{eq:poinc:31}
    g(\lambda w)=\tilde{R}_{\lambda,r}(g(w))=\lambda \tilde{S}_r(g(w))=\lambda\int_0^{g(w)}(1-u)^{n}(1-\frac{u}{r})(1+\frac{u}{r^r})du,\quad r>1.
    \end{equation}
    In this case
    \[
    2S_r(r)\geq \tilde{S}_r(r)\geq S_r(r).
    \]
    Thus we may again set $0<\lambda=r/\tilde{S}_r(r)<1$, so that $\tilde{R}_{\lambda,r}(r)=r$, and as before we can find $\lambda_4$ such that
    \[
    \tilde{R}_{\lambda_5,r}^{\circ k}(r)<1<\tilde{R}_{\lambda_5,r}^{\circ (k-1)}(r)
    \]
    for all $\lambda_5\in(\lambda_4,\lambda_4+\epsilon)$ for some small enough $\epsilon>0$. Again we set
    \[
    \lambda_6=\inf_{t\in(0,r]}\frac{t}{\tilde{S}_r(t)}>0.
    \]
    Arguing as for $S_r$, we obtain that
    \[
    \lambda_6\in(\frac{b_1}{r^n},\frac{b_2}{r^n}).
    \]
    The proofs of $\tilde{R}_{\lambda_6,r}(r)\geq v_1r>1$ and $\tilde{R}_{\lambda_6,r}^{\circ j}(r)>u_jr>1$ are the same as in Example \ref{ex:poinc:6}, and so is the proof that $\tilde{R}_{\lambda_5,r}^{\circ (k-1)}(r)$ is not a value of $g_0$. This gives the example for polynomials of odd degree.
\end{example}

\section{Main lemmas for FDEs}
\label{sec:nev}

\subsection{Classical Nevanlinna theory}
\label{subsec:nev:1}

In the next two sections we study the singularity structure of the analytic continuation of the solutions of (\ref{eq:intro:7}) and (\ref{eq:intro:8}). In the generic case, we show that these solutions cannot be meromorphic, which allows us to prove the existence of natural boundaries or branch points. For this we use basic results of classical Nevanlinna theory. We recall the definitions and some important theorems without proof; see \cite{hayman1964meromorphic,laine93,golʹdberg2008value}.

\begin{definition}
\label{def:nev:1}
    Let $f$ be meromorphic function, and $n(r,f)$ be the number of poles of $f$ counted with multiplicities inside $\{t\in\mathbb{C}:|t|\leq r\}$. We define
    \begin{align*}
    &\mbox{Nevanlinna counting function:} & &N(r,f)=\int_0^r\frac{n(s,f)-n(0,f)}{s}\mathrm{d} s+n(0,f)\ln r,\\
    &\mbox{Proximity function:} & &m(r,f)=\frac{1}{2\pi}\int_0^{2\pi}\ln^+|f(r\e^{i\theta})|\mathrm{d}\theta.
    \end{align*}
    We define Nevanlinna characteristic by
    \[
    T(r,f)=m(r,f)+N(r,f).
    \]
We denote $U(r,a,f):=U(r,1/(f-a))$ for $a\in\mathbb{C}$, where $U$ can be either $m,~N$ or $T$ above; $U(r,f)$ is sometimes written as $U(r,\infty,f)$.
\end{definition}

\begin{theorem}[Nevanlinna's First Main Theorem]
\label{thm:nev:1}
    Let $a\in\mathbb{C}$, then
    \[T(r,f)=T(r,a,f)+O(1)\quad
    \mbox{\ as\ \  }r\to\infty.
    \]
\end{theorem}

\begin{theorem}
\label{thm:nev:2}
Let $U$ be $N$, $m$ or $T$ in Definition \ref{def:nev:1}, and let $f,g$ be meromorphic.  Then $U(r,fg)\leq U(r,f)+U(r,g)+O(1),~ U(r,f+g)\leq U(r,f)+U(r,g)+O(1)$.
\end{theorem}

\begin{definition}
\label{def:nev:2}
For functions $P_1(r),P_2(r)$, if there exist $\omega_1,\omega_2>0$ such that $\omega_1 P_1\leq P_2\leq\omega_2 P_1$ for all large $r$, we write $P_1\asymp P_2$.

Most results in Nevanlinna theory hold outside a set of finite Lebesgue measure. We also use the notation $P_1\asymp P_2$ when the inequalities hold only outside such a set.
\end{definition}

\begin{theorem}
\label{thm:nev:3}
    Let $f$ be meromorphic. Then
    \begin{itemize}
        \item $f$ is a constant if and only if $T(r,f)=O(1)$,
        \item $f$ is a non-constant rational function if and only if $T(r,f)\asymp \ln r$,
        \item $f$ is a transcendental function if and only if $\ln r=o(T(r,f))$.
    \end{itemize}
\end{theorem}

\begin{theorem}
\label{thm:nev:4}
    A meromorphic function $f$ is transcendental if and only if $A(r,f):=T(r,f)/\ln r\to\infty$ as $r\to\infty$. In this case, $A(r,f)$ is increasing for all large $r$.
\end{theorem}

\begin{definition}
\label{def:nev:3}
Let $G(r)$ be a non-negative function such that $G(r)=o(T(r,f))$ as $r\to\infty$ outside a set of finite Lebesgue measure\footnote{also known as linear measure}. We denote $G(r)=S(r,f)$. Any such function $G$ is called a small function with respect to $f$ in the Nevanlinna sense.  We also call a meromorphic function $g$ small with respect to $f$ if $T(r,g)$ is small with respect to $f$.
\end{definition}

The following theorem is a strong form of Lemma on the Logarithmic derivative (LOTLD) \cite[Ch.2]{golʹdberg2008value}.
\begin{theorem}
\label{thm:nev:5}
Given $0<\alpha<1$ and a meromorphic function $f$, then
\begin{equation}
\label{eq:nev:1}
m\left(r,\frac{f'}{f}\right)\leq\frac{\ln^+T(R,f)}{\alpha}+\max\left\{2,\frac{1}{\alpha}\right\}\left(\ln^+\frac{1}{R-r}+\ln R\right)+O(\ln R).
\end{equation}
for all $1<r<R<\infty$. Moreover,
\[
\sup\left\{m\left(r,\frac{f^{(n)}}{f^{(n-1)}}\right),m\left(r,\frac{f^{(n)}}{f}\right)\right\}\leq O\left(\sum_{j=1}^{n}\ln^{\circ j+}T(R_{j},f)+\ln R_n\right)
\]
for all large $r<R_1<\cdots<R_n$, where $\ln^{\circ j+}=\ln^+\circ\cdots\circ\ln^+~(j\text{-times})$.
\end{theorem}

The LOTLD is usually stated as follows.
\begin{theorem}[LOTLD]
\label{thm:nev:6}
Let $f$ be a non-constant meromorphic function. Then
\[
m(r,f'/f)=S(r,f).
\]
Moreover, $m(r,f^{(k)}/f)=S(r,f)$ for all $k\in\mathbb{N}$. Hence, $T(r,f^{(n)})=O(T(r,f))$.
\end{theorem}

\begin{theorem}[Valiron-Mohon'ko]
\label{thm:nev:7}
Let $R(f)=P(f)/Q(f)$ where $P(f)=a_0+a_1f+\cdots+a_pf^p$ and $Q(f)=b_0+b_1f+\cdots+b_qf^q$ such that $\gcd(P,Q)=1$. If $T(r,a_i)=S(r,f)$ and $T(r,b_j)=S(r,f)$ for $i=0,\dots,p$ and $j=0,\dots,q$, then $T(r,R(f))=\max\{p,q\}T(r,f)+S(r,f)$.
\end{theorem}

\subsection{Existence of Laurent series}
\label{subsec:nev:2}

In Sections \ref{sec:fde} and \ref{sec:qde} we discuss solutions with a simple pole at $z=0$ of the FDEs (\ref{eq:intro:7}) and (\ref{eq:intro:8}). In other words, we look for solutions of the form
\begin{equation}
    \label{eq:nev:2}
y=\frac{\theta_{-1}}{z}+\sum_{j=0}^\infty \theta_jz^j,
\end{equation}
where the $\theta_j$ are constants. We need to show that the analytic part $\sum_{j=0}^\infty \theta_jz^j$ of the formal series converges. The proof is standard and elementary, so we present it in Appendix~\ref{sec:appA}: we construct a majorising series $\hat{y}$, analytic about $z=0$, whose coefficients $\hat{\theta}_j$ bound the coefficients $\theta_j$ (see \cite[Ch.2]{wasow2018asymptotic} and \cite[§2.4]{hille1997ordinary}).

\begin{theorem}
\label{thm:nev:8}
Let $y$, given by (\ref{eq:nev:2}), be a formal Laurent series solution of either
\begin{equation}
\label{eq:nev:3}
y'(z)+a(z)y^2(z)+b(z)y(z)+d(z)y(z^2h(z))+f(z)=0,
\end{equation}
or
\begin{equation}
\label{eq:nev:4}
{\begin{aligned}
&y'(z)+a(z)y^2(z)+b(z)y(z)+c(z)y(z\check{h}(z))\\&+d(z)y^2(z\check{h}(z))+e(z)y(z)y(z\check{h}(z))+f(z)=0.
\end{aligned}}
\end{equation}
Here $a,b,c,d,e,f,h$ and $\check{h}$ are analytic at $0$, with $h(0)\neq0$, $\check{h}(0)\neq0$ and $|\check{h}(0)|<1$; when $\check{h}(z)\equiv q$ is constant, so that $y(z\check{h}(z))=y(qz)$, we allow $|q|\leq1$. We write $a(z)=\sum_{j=0}^\infty a_jz^j$ for the Taylor series of $a$, and similarly for $b,c,d,e,f,h$ and $\check{h}$.
If $a_0,d_0,e_0$ are such that $\theta_{-1}\neq0$ and $\Srec(n)\neq0$ for all $n\geq-1$, then the analytic part of (\ref{eq:nev:2}) converges in a small disc around $z=0$.
\end{theorem}

Here, $\theta_{-1}=(h_0-d_0)/(a_0h_0)$ in (\ref{eq:nev:2}) for equation (\ref{eq:nev:3}) and $\theta_{-1}=\check{h}_0^2/(a_0\Check{h}_0^2+e_0\Check{h}_0+d_0)$ in (\ref{eq:nev:2}) for (\ref{eq:nev:4}). Hence, we need the additional assumptions $a_0(h_0-d_0)\neq0$ for (\ref{eq:nev:3}) and $a_0\Check{h}_0^2+e_0\Check{h}_0+d_0\neq0$ for (\ref{eq:nev:4}). In fact, if $h_0=d_0$, (\ref{eq:nev:2}) still exists and $a_0$ is zero\footnote{this is the resonance condition, where $\theta_{-1}$ is arbitrary.}.
Now we assume that $\theta_{-1}\neq0$. To show the existence, we observe the recurrence relation:
\begin{equation}
\label{eq:nev:5}
\Srec(n)\theta_{n+1}=\Trec(\theta_{-1},\theta_0,\theta_1,\dots,\theta_{n})
\end{equation}
where $\Srec$ and $\Trec$ are polynomials and $\Srec(n)$ depends only on $a_0,~d_0,~e_0$. The zeros of $\Srec(n)$ are the resonances. When $\Srec(n)\neq0$, the coefficients $\theta_0,\theta_1,\dots$ are determined uniquely by the recurrence relation (\ref{eq:nev:5}). When $\Srec(n)=0$, we obtain the resonance condition $\Trec(\theta_{-1},\theta_0,\theta_1,\dots,\theta_{n})=0$. The following corollary treats the case $\Srec(n)=0$.

\begin{cor}
\label{cor:nev:1}
If the condition $\Srec(n)\neq0$ in Theorem \ref{thm:nev:8} is relaxed, then there are at most finitely many $n$ such that $\Srec(n)=0$. In this case, the conclusion of Theorem \ref{thm:nev:8} remains true, subject to the resonance conditions.
\end{cor}

\begin{remark}
    The existence of an analytic solution at $z=0$ of (\ref{eq:nev:3}) and (\ref{eq:nev:4}) is an immediate consequence of Oberg's result \cite[Th.2.1]{oberg1971local} on FDEs near an attracting fixed point; the authors of \cite{marshall2002natural,marshall2004natural} use this result to prove the existence of an analytic solution at $z=0$ of a pantograph equation.
\end{remark}


\subsection{Non-meromorphic analytic continuation}
\label{subsec:nev:3}

\begin{theorem}
\label{thm:nev:9}
Let $y$ be a transcendental meromorphic function. Let $L$ be the operator defined by
\[
L[y]=\sum_{j=0}^qa_jy^j(y')^{j_1}\cdots(y^{n_j})^{j_{n_j}},
\]
where the $a_j$ are meromorphic functions such that $T(r,a_j)=o(T^*(r,y))$. Here, $T^*(r,y)=T(r,y)/(\ln\ln r)^{1+\epsilon}$ for some $\epsilon>0$. We may rewrite $L$ as $L[y]=\sum_{i=0}^\tau y^iP_i(y'/y,\dots,y^{(n_i)}/y)$, where $P_i$ is a polynomial in $y'/y,\dots,y^{(n_i)}/y$ whose coefficients $a_{il}$ satisfy $T(r,a_{il})=o(T^*(r,y))$. Let $g$ be an entire function that is neither of the form $az+b$ with $b\neq0$ nor of the form $\gamma z$ with $|\gamma|\leq1$. If $y$ satisfies an equation
\[
L[y]=P(y\circ g (z)),
\]
where $P(t)=\sum_{i=0}^\sigma c_it^i$ with $T(r,c_i)=o(T^*(r,y))$ and $c_\sigma=1$, then $g$ must be a polynomial. Furthermore, if $N(r,L[y])\leq \tau N(r,y)+o(T^*(r,y))$, then the following hold:
\begin{enumerate}
    \item when $g(z)=\gamma z,~|\gamma|>1$, then $\tau\geq\sigma$.
    \item when $\deg g=\mu\geq 2$, then $\tau>\sigma\mu$.
\end{enumerate}
If the conditions $o(T^*(r,y))$ are replaced by $o(T(r,y))$, then the above statements remain true, except that (2) must be replaced by:
    \begin{enumerate}
        \item[(2$'$)] when $\deg g=\mu\geq 2$, then either $\tau>\sigma\mu$, or $\tau=\sigma\mu$ and $y$ has order $\rho=\limsup_{r\to\infty}\ln T(r,y)/\ln r=0$.
    \end{enumerate}
\end{theorem}
The condition $c_\sigma=1$ only guarantees that $c_\sigma\not\equiv0$. If $c_\sigma$ is a meromorphic function satisfying the same growth condition as the other coefficients, we can divide both sides of $L[y]=P(y\circ g(z))$ by $c_\sigma$.

\begin{proof}
First, we show that $g$ is not transcendental. Assume the contrary. By Theorem \ref{thm:nev:5},
\begin{equation}
\label{eq:nev:6}
T(r,L[y])\leq O( T(r,y))+O\left(\sum_{j=1}^v\ln^{\circ j+}T(R_j,y)+\ln R_v\right)\leq O(T(R_v,y)),
\end{equation}
where we choose $R_v=r+\frac{1}{4T(r,y\circ g)}$. On the other hand, $T(r,L[y])=T(r,P(y\circ g))=\sigma T(r,y\circ g)+o(T(r,y))$ by Theorem \ref{thm:nev:7}. Together with (\ref{eq:nev:6}), this gives
\[
T(r,y\circ g)\leq O\left(T\left(r+\frac{1}{T(r,y\circ g)},y\right)\right),
\]
which contradicts (\ref{eq:appB:8}) in Appendix \ref{sec:appB}. Hence $g$ is a polynomial.

Let $g=b_\mu z^\mu+\cdots+b_0$ with $\mu\geq2$ and $b_\mu\neq0$. Taking the characteristic function of both sides of $P(y\circ g)=L[y]$, and using $N(r,L[y])\leq \tau N(r,y)+o(T(r,y))$ and Theorem \ref{thm:nev:5}, we obtain
\begin{equation}
\label{eq:nev:7}
    \sigma T(r,y\circ g)\leq \tau T(r,y)+O(\ln^+T(R_v,y))+O(\ln R_v)+o(T(r,y))
\end{equation}
for any $R_v>r$. By the left-hand inequality in (\ref{eq:appB:4}), $T(r,y\circ g)\geq T(d_1r^\mu,y)-K_1$. We choose $R_v=d_1r^\mu$; since $\ln^+T(R,y)=o(T(R,y))$ and $\ln r=o(T(r,y))$, (\ref{eq:nev:7}) gives
\begin{equation}
\label{eq:nev:8}
    (\sigma-o(1))\,T(d_1r^\mu,y)\leq \tau T(r,y)+o(T(r,y)).
\end{equation}
On the other hand, $A(r,y)=T(r,y)/\ln r$ is increasing for all large $r$ by Theorem \ref{thm:nev:4}, so
\begin{equation}
\label{eq:nev:9}
T(d_1r^\mu,y)\geq\frac{\mu\ln r+\ln d_1}{\ln r}\,T(r,y)=(\mu+o(1))\,T(r,y).
\end{equation}
If $\tau\leq\sigma\mu-1$, then (\ref{eq:nev:8}) and (\ref{eq:nev:9}) give $(\sigma\mu-o(1))T(r,y)\leq(\sigma\mu-1+o(1))T(r,y)$, a contradiction. Hence $\tau\geq\sigma\mu$.

Suppose now that $\tau=\sigma\mu$. Then, by (\ref{eq:nev:8}), for every $\delta>0$ we have
\begin{equation}
\label{eq:nev:10}
    T(d_1r^\mu,y)\leq \mu(1+\delta)\,T(r,y)\qquad\text{for all }r\geq R_0,
\end{equation}
for some $R_0=R_0(\delta)$. Let $\rho=\limsup_{r\to\infty}\ln T(r,y)/\ln r$ be the order of $y$. Taking logarithms in (\ref{eq:nev:10}), we obtain
\[
\frac{\ln T(d_1r^{\mu},y)}{\ln(d_1r^\mu)}\leq\frac{\ln T(r,y)+O(1)}{\mu \ln r+\ln d_1},
\]
which yields $\rho\leq\rho/\mu$. Thus $\rho=0$ or $\rho=+\infty$.

We show that $\rho=+\infty$ is impossible. The argument does not depend on $d_1$, so for convenience we take $d_1=1$; in general, $r^\mu$ is replaced by $d_1r^\mu$ throughout. Fix $r_0\geq R_0$ such that $r_0^{\mu-1}\geq2\mu(1+\delta)$. By (\ref{eq:nev:10}), for every $r\geq r_0$,
\[
\frac{T(r^\mu,y)}{r^\mu}\leq\frac{\mu(1+\delta)}{r^{\mu-1}}\cdot\frac{T(r,y)}{r}\leq\frac12\cdot\frac{T(r,y)}{r}.
\]
Hence, by induction, for every $s\in[r_0,r_0^\mu]$ and every $n\geq0$,
\begin{equation}
\label{eq:nev:11}
\frac{T(s^{\mu^n},y)}{s^{\mu^n}}\leq2^{-n}\,\frac{T(s,y)}{s}\leq2^{-n}\,\frac{T(r_0^\mu,y)}{r_0}.
\end{equation}
Now suppose that $\rho=+\infty$. Then for every $N>0$ there exist arbitrarily large $r$ with $T(r,y)/r\geq N$. Every $r\geq r_0$ can be written as $r=r_0^{\mu^{p}}$ with $p\geq0$, and then $r=s^{\mu^n}$ with $n=\lfloor p\rfloor$ and $s=r_0^{\mu^{p-\lfloor p\rfloor}}\in[r_0,r_0^\mu)$. By (\ref{eq:nev:11}),
\[
N\leq\frac{T(r,y)}{r}\leq2^{-\lfloor p\rfloor}\,\frac{T(r_0^\mu,y)}{r_0},
\]
and the right-hand side tends to $0$ as $r\to\infty$, a contradiction. Hence $\rho\neq+\infty$, so $\rho=0$. Thus either $\tau>\sigma\mu$, or $\tau=\sigma\mu$ and $\rho=0$; this proves (2$'$).

Now we prove the case $\tau=\sigma\mu$ and $\rho=0$.
Here we use the stronger condition $o(T^*(r,y))$ on the coefficients. Recall the following strong form of the LOTLD, without an exceptional set \cite[Exercise 1, p.91]{golʹdberg2008value}: if $f$ is a meromorphic function of order less than $1$, then $m(r,f^{(n)}/f)=o(1)$ as $r\to\infty$. Since $\rho=0$, then
\begin{equation}
\label{eq:nev:12}
    m\left(r,\frac{y^{(n)}}{y}\right)=o(1).
\end{equation}
By (\ref{eq:nev:12}), Theorem \ref{thm:nev:7} and the hypotheses on the coefficients,
\begin{align*}
\sigma T(r,y\circ g)&=T(r,P(y\circ g))+o(T^*(r,y))\\
&=T(r,L[y])+o(T^*(r,y))\\
&\leq \tau T(r,y)+o(T^*(r,y)).
\end{align*}
Since $\tau=\sigma\mu$, the left-hand inequality in (\ref{eq:appB:4}) gives
\begin{equation}
\label{eq:nev:13}
    T(d_1r^{\mu},y)\leq \mu T(r,y)+o(T^*(r,y)).
\end{equation}
Dividing by $\ln(d_1r^\mu)=\mu\ln r+\ln d_1$, we obtain, for all large $r$,
\begin{equation}
\label{eq:nev:14}
    A(d_1r^{\mu},y)\leq A(r,y)\left(1+\frac{-\ln d_1}{\mu\ln r+\ln d_1}\right)\left(1+\frac{1}{(\ln\ln r)^{1+\epsilon}}\right).
\end{equation}
Let $r_0$ be large with $d_1r_0^{\mu-1}>1$, and let $r_{i+1}=d_1r_i^{\mu}$ for $i\geq0$. Then $r_i$ increases to $\infty$, and $\ln r_i+\beta=\mu^i(\ln r_0+\beta)$ with $\beta=\ln d_1/(\mu-1)$, so $\ln r_i$ grows geometrically and $\ln\ln r_i=i\ln\mu+O(1)$. Applying (\ref{eq:nev:14}) with $r=r_0,\dots,r_{k-1}$, we obtain
\begin{equation}
\label{eq:nev:15}
    A(r_k,y)\leq A(r_0,y)\prod_{i=0}^{k-1}\left(1+\frac{-\ln d_1}{\mu\ln r_i+\ln d_1}\right)\left(1+\frac{1}{(\ln\ln r_i)^{1+\epsilon}}\right).
\end{equation}
The series
\[
\sum_{i=0}^{\infty}\frac{|\ln d_1|}{\mu\ln r_i+\ln d_1}\quad\text{and}\quad\sum_{i=0}^\infty\frac{1}{(\ln\ln r_i)^{1+\epsilon}}
\]
both converge, so the infinite product converges, and $A(r_k,y)=O(A(r_0,y))$ for all $k$. Since $r_k\to\infty$, this contradicts $A(r,y)\to\infty$ (Theorem \ref{thm:nev:4}). Hence $\tau>\sigma\mu$, which proves (2).

Finally, we discuss the case $g(z)=\gamma z$ with $|\gamma|>1$. Suppose that $\tau+1\leq \sigma$. Arguing as above, we obtain
\begin{align*}
&T(r,P(y\circ g))=\sigma T(r,y(\gamma z))+o(T(r,y))\geq (\tau+1)T(r|\gamma|,y)+o(T(r,y)),\\
&T(r,L[y])\leq \tau T(r,y)+O(\ln^+T(R_v,y))+o(T(r,y)).
\end{align*}
Choosing $R_v=|\gamma|r$ gives $(\tau+1-o(1))T(r|\gamma|,y)\leq \tau T(r,y)$, a contradiction. Hence $\tau\geq\sigma$, which proves (1).
\end{proof}

\begin{remark}
\label{rem:nev:2}
The above remains valid if $T^*(r,y)=u(r)$, where $u$ is a positive continuous function such that $\sum_{i=0}^\infty u(r_i)/T(r_i,y)<\infty$, where $r_{i+1}=d_1r_i^\mu$ for some large $r_0$.
\end{remark}


Now we generalize the results in \cite{marshall2002natural,marshall2005holomorphic} to some nonlinear FDEs. Theorem \ref{thm:nev:10} and Corollary \ref{cor:nev:2} are non-linear analogues of the results in \cite[Lemmas 4.25 and 4.34, Corollaries 4.26 and 3.35, and Theorem 4.31]{marshall2005holomorphic}.
We first introduce some definitions and facts.

\begin{definition}
\label{def:nev:4}
Let $R$ be an entire function or a rational function. A sequence of functions $\{R^{\circ n}\}$ is said to be normal at $w_0$ if there exists a neighborhood $D$ of $w_0$ such that there is a subsequence of $\{R^{\circ n}(z)\}$ converging uniformly in $D$; here it also includes the convergence to $\infty$ using the spherical distance, however, we can map this point to a finite point and the spherical distance can be replaced locally by the usual complex modulus.

The Julia set of $R$ is the set $J(R)$, which contains all the points $w\in\mathbb{C}$ such that the sequence $(R^{\circ n})$ is not normal in any neighborhood of $w$. When $R$ is a rational function, the Julia set also includes the possibility of the point at infinity.
\end{definition}

We state some facts about the Julia set $J(R)$ without proof, see \cite[Ch.3]{carleson2013complex}, \cite[§3-5]{milnor2011dynamics}, \cite[Ch.3-4]{beardon2000iteration}.
\begin{enumerate}
    \item \underline{Fact 1:} $J(R)$ is closed (non-empty if $R$ is not a linear function).
    \item \underline{Fact 2:} The set of points $w\in J(R)$ such that the orbit $\{R^{\circ n}(w),~n\geq0\}$ is dense on $J(R)$ has full measure on $J(R)$.
    \item \underline{Fact 3:} Given any $w\in J(R)$ and any neighborhood $U$ of $w$, then $\cup_{n\geq0} R^{\circ n}(U)$ omits at most one point in $\mathbb{C}$. The omitted value is known as an exceptional value, which always lies in the Fatou set $F(R)$; here $F(R):=\mathbb{\hat{C}}\setminus J(R)$ if $R$ is a rational function, and $F(R):=\mathbb{C}\setminus J(R)$ if $R$ is a transcendental entire function.
    \item \underline{Fact 4:} Given $w\in J(R)$ and any point $w_0$ that is not an exceptional value, then $w$ is a limit point of $\mathbf{L}:=\{(R^{-1})^{\circ n}(w_0),n\geq0\}$ (of any branches).
    \item \underline{Fact 5:} Given any $w\in J(R)$, then $\mathbf{L}$ is dense in $J(R)$. Outside a set of measure zero, any $w\in J(R)$ we have $\{R^{\circ n}(w),n\geq0\}$ is dense in $J(R)$.
\end{enumerate}

\begin{definition}
\label{def:nev:5}
Let $f$ be a function defined on a domain $D$, except on some countable set in the interior of $D$. If $f$ is analytic at a point $z_0\in \mathbb{C}$, then $z_0$ is said to be a regular point of $f$. If none of the points on $\partial D$, the boundary of $D$, is regular, then $\partial D$ is called the natural boundary of $f$. For a multivalued $f$, we say that $\partial D$ is a natural boundary of $f$ if no branch of $f$ can be continued analytically, along a path in $D$, to a point of $\partial D$.
\end{definition}

\begin{theorem}
\label{thm:nev:10}
Let $g$ be either a polynomial of degree $\mu\geq2$ or a transcendental entire function, with an attracting fixed point $z_0$. Assume further that each $z\in J(g)$ has a neighbourhood $U_z$ such that $\cup_{n\geq0}g^n(U_z)$ is either $\mathbb{C}$ or $\mathbb{C}\setminus\{z_0\}$.
Suppose that $y$ is a meromorphic solution of $L[y]=y(g(z))$ in a neighbourhood of $z_0$, where $L$ is as in Theorem \ref{thm:nev:9}, and that $\mu\geq \tau$ when $g$ is a polynomial. If $y$ is transcendental, then the analytic continuation of $y$ on any Riemann sheet cannot be extended to any point of $\partial A_{z_0}^0$.
\end{theorem}

\begin{cor}
\label{cor:nev:2}
Let $g$, $z_0$, $y$ and $L$ be as in Theorem \ref{thm:nev:10}, and assume in addition that each $z\in J(g)$ has a neighbourhood $U_z$ such that $\cup_{n\geq0}g^n(U_z)=\mathbb{C}\setminus\{\acute{z}\}$ for some $\acute{z}\neq z_0$. Then the following hold.
\begin{enumerate}
    \item If $y$ has a pole at $z_0$, then the conclusion of Theorem \ref{thm:nev:10} still holds.
    \item If $y$ is analytic at $z_0$, then either the conclusion of Theorem \ref{thm:nev:10} holds, or $y$ is holomorphic on $\mathbb{C}\setminus\{\acute{z}\}$ and $\acute{z}$ is an essential singularity of $y$.
    \item In both cases, if the conclusion of Theorem \ref{thm:nev:10} holds and $y$ is transcendental, then $y$ continues to a Riemann surface lying over $A_0^0$ minus a discrete set, and $y$ cannot be extended to $\partial A_{z_0}^0$ on any Riemann sheet.
\end{enumerate}
\end{cor}

\begin{proof}[Proof of Theorem \ref{thm:nev:10}]
We can assume that $z_0=0$. Theorem \ref{thm:nev:9} shows that the Laurent series of $y$ at $0$ has a radius of convergence $0<R<\infty$. So there is a point $c_0\in C(0,R)$ such that $y$ is singular about $c_0$. If $y$ can be continued such that $y$ is analytic at a point $c\in J(g)$, then there exists a small disk $D(c,\epsilon)$ on which $y$ is analytic. From Fact 4 (below Definition \ref{def:nev:4}), we find that $y$ is analytic on $\mathbb{C}-\{z_0\}$. This shows that $y$ is meromorphic on $\mathbb{C}$, which contradicts Theorem \ref{thm:nev:9}. Hence, $y$ cannot be developed at any points on $J(g)$. So the point $c_0$ must be in $\overline{A_0^0}$. If $c_0\in\partial A_0^0$, we continue $y$ through $A_0^0$. If there are no other singular points, then $y$ is single-valued in $A_0^0$. By the above argument, $\partial A_0^0$ is a natural boundary of $y$. If there exists a point $c_1\in A_0^0$ (including $c_0=c_1$), then $y$ is singular on $\{g^{-\circ n}(c_1),~n\geq0\}$ which accumulates to every point on $\partial A_0^0$. This shows that $\partial A_0^0$ is the natural boundary of $y$. If none of the points $\{g^{-\circ n}(c_1),~n\geq0\}$ are branched, then $y$ is single-valued. Now we consider when one of the points is branched, we may assume that $c_1$ is branched. Then at $c_1$, $y$ yields the sheets $y_i$, where $y_0$ is the predecessor sheet and the rest are the subsequent successor sheets (see Definition \ref{def:poinc:1}). We also have $L[y_i]=y_0(g(z))$, which shows that $y_i$ are also branched at $\{g^{-\circ n}(c_1),~n\geq0\}$ which shows that $y_i$ admits $\partial A_0^0$ as the natural boundary. From there, the succesor sheets $y_{ij}$ of each $y_i$ also satisfy $L[y_{ij}]=y_i$, which shows that $g^{-\circ n}(c_1),~n\geq1$ are branch points of $y_{ij}$. Hence, $\partial A_0^0$ is the natural boundary of $y_{ij}$. We repeat this argument, which asserts that the succesor sheets occurred from the branch points $g^{-\circ n}(c_1)$ admit $\partial A_0^0$ as the natural boundary. On those sheets, say $y_{ij}$ could be branched on other points not on $\{g^{-\circ n}(c_1),~n\geq0\}$, we may assume that the point is $w_0$. The successor sheets developed from the point $w_0$, denoted by $\tilde{y}_{ij}$, which satisfies $L[\tilde{y}_{ij}]=y_i(g(z))$. This shows that the set $\{g^{-\circ n}(c_1),~n\geq N\}$ for large enough $N$, are also branch points of $y_{ij}'$. Hence, $\tilde{y}_{ij}$ has the natural boundary $\partial A_0^0$.

To show the last statement, we see that when we solve $L[\tilde{y}_{i}]=y(g(z))$, there may be some undefined discrete points. When we continue $y$ from a neighborhood of $z_0$, we use $L[\tilde{y}_{i}]=y(g(z))$ to solve $y$; this neighborhood is in $A_0^0$ as Theorem \ref{thm:nev:10} shows that we cannot develop beyond $\partial A_0^0$. In this procedure, there could be a singularity or undefined point. For instance, if $L[y]=y'(z)+ay^2(z)$ and $y$ has a pole at $z_1$, then there could be a logarithmic branch point at $g^{-1}(z_1)$. These undefined points or singularities form a discrete set and could occur on any stage of development of $y$ on any Riemann sheets. Assume that on $y_i$, the set of those points is $\tilde{w}_\beta$, where $\beta$ is an index. Then $y_i$ is not defined on $\{g^{-\circ n}(\tilde{w}_\beta),~n\geq0\}$ for any index $\beta$. This undefined set consists of discrete points and will be inherited through the subsequent successor sheets of $y_i$. For other sheets, the argument is similar.
\end{proof}

\begin{remark}
\label{rem:nev:4}
Theorem \ref{thm:nev:10} does not state that $\partial A_{z_0}^0$ is a natural boundary of $y$ due to complexity of $L$, unlike when $L$ is linear, which has been shown in \cite{marshall2002natural}. However, there is indeed a natural boundary of $y$ on any sheet which is contained in $\overline{A_{z_0}^0}$. In section \ref{sec:fde}, we will provide more concrete results when $L[y]=y'+P(y(g(z)),y(z))$ where $P$ is a degree $2$ polynomial.
\end{remark}

\begin{proof}[Proof of Corollary \ref{cor:nev:2}]
When $y$ has a pole at $z_0$, we can modify the proof of Theorem \ref{thm:nev:10}. Again, we assume $z_0=0$. Since $0$ is not an exceptional value, then the argument in Theorem \ref{thm:nev:10} still applies: if there is a point $c\in J(g)$ that is a regular point of $y$ i.e. $y$ is holomorphic on $U_c$, then $y$ is analytic on $\cup_{n\geq0}g^n(U_c)=\mathbb{C}-\{\Acute{z}\}$. Moreover, there exists $z\in U_c$ such that $g^{n_0}(z)=0$ and so $y$ is regular at $0$ contradicting the fact that $0$ is a pole of $y$. From there, we argue similarly as in Theorem \ref{thm:nev:10}.

Now we consider the case where $y$ is analytic at $z_0=0$. We assume that the conclusion of Theorem \ref{thm:nev:10} does not hold. Now the function $y$ in the proof of Theorem \ref{thm:nev:10} is holomorphic everywhere except at $\acute{z}$. Since $y$ is single-valued on a punctured disk centered on $\acute{z}$, then $\acute{z}$ is an (isolated) essential singularity.
\end{proof}

\section{Singularity structure of special solutions of the FDE (\ref{eq:intro:7})}
\label{sec:fde}

In this section, we discuss solutions of (\ref{eq:intro:7}) that have poles at $0$. There also exist solutions that are analytic at $0$, which follow from \cite[Th.2.1]{oberg1971local}. However, our attention will be given to solutions with poles at $0$, and we will remark on the analogous results in the analytic cases; when $b=-1,a=0$, the authors in \cite{marshall2005holomorphic} used the Hadamard gap theorem \cite[Ch.6,§6]{segal2007nine} to show that the analytic solution at $0$ has the radius of convergence $1$ and $C(0,1)$ is the natural boundary. Theorem \ref{thm:nev:8} and Corollary \ref{cor:nev:1} show that there exists a solution $y$ of (\ref{eq:intro:7}) with a simple pole at $0$. In this paper we only consider such solutions, which give the simplest description of the singularity structure; poles of higher order are discussed in Section \ref{sec:disc}. Thus
\begin{equation}
\label{eq:fde:1}
y=\frac{\theta_{-1}}{z}+\sum_{i=0}^\infty\theta_iz^i,\quad \theta_{-1}\neq0.
\end{equation}
By Theorem \ref{thm:nev:10}, if (\ref{eq:fde:1}) is transcendental, then its analytic continuation cannot be developed to any point of $\partial A_0^0=J(z^2)=C(0,1)$; for the map $z^2$, the points $0$ and $\infty$ are the exceptional values.
The main result is the following theorem.

\begin{theorem}
\label{thm:fde:1}
Let (\ref{eq:fde:1}) be the solution of (\ref{eq:intro:7}). Analytic continuation of (\ref{eq:intro:7}) leads to the following results:
\begin{enumerate}
    \item if $b=1$, then $a=0$, and either $y=\theta_{-1}/z$, or $y$ is a solution of (\ref{eq:intro:7}) with $\theta_0\neq0$, which is a single-valued transcendental function defined on $D(0,1)$.
    \item if $b\neq1$ and $b\neq(N+2)/2$ for all $N\in\mathbb{N}$, then $y=(1-b)/(az)$.
    \item if $b=(N+2)/2$ for some $N\in\mathbb{N}$, then either $y=(1-b)/(az)$, or $y=(1-b)/(az)-bz^{N}/a$, or $y=(1-b)/(az)+\sum_{j=N}^\infty \theta_jz^j$ with $\theta_N\neq0,~-b/a$, which is transcendental.
    \item Each transcendental solution can be continued meromorphically into the unit disc, as a possibly multivalued function, with the exception of some singularities, and each sheet has the natural boundary $C(0,1)$.
    \item Let $\tilde{r}$ be the radius of convergence of the analytic part of (\ref{eq:fde:1}), for transcendental $y$. If $\tilde{r}=1$, then $y$ is a meromorphic function on $D(0,1)$. If $\tilde{r}<1$, then $y$ cannot be a single-valued meromorphic function on $D(0,1)$; that is, $y$ is multivalued.
\end{enumerate}
\end{theorem}

\begin{proof}
By substituting (\ref{eq:fde:1}) into (\ref{eq:intro:7}), we obtain the recurrence relation. If $b=1$, we find easily that $a=0$. In this case, $\theta_{-1}$ is a parameter, and (\ref{eq:intro:7}) becomes $y'+y(z^2)=0$. It is easy to see that $y$ only has a pole at $z=0$, otherwise, if $z_0$ is a pole of $y$, then the sequence $\{z_0^{1/2^n},~n\geq0\}$ is a sequence of poles of $y$ accumulating to every point on $C(0,1)$. Hence, $y$ admits only a pole at $z=0$. The recurrence relation is $\theta_{2k+2}=0,~(2k+1)\theta_{2k+1}+\theta_k=0,~k\geq0$. Since $y$ has only a pole at $z=0$, then when $y$ is rational, $y=P(z)/z$ for some polynomial $P$. In this case, $\theta_k=0$ for all large $k$. However, from the recurrence relation, we see that this occurs only if $\theta_k=0$ for $k\geq0$ which corresponds to $y=\theta_{-1}/z$. For this not to occur, we must have $\theta_0\neq0$, and $y$ is transcendental. When $y$ is transcendental, we can see from $y'(z)=-y(z^2)$ that $y$ admits only a pole at $z=0$; otherwise, $z_0$ is a pole and we obtain a sequence $\{z_0^{2^n}\}$ of poles accumulating to $0$. Thus, from solving $y'(z)=-y(z^2)$, we find that $y$ is defined on $D(0,1)$ with only a pole at $z=0$. Hence, Theorem \ref{thm:nev:10} shows that $C(0,1)$ is a natural boundary of $y$.

Now we assume that $b\neq1$. It can be checked that $y=(1-b)/(az)$ is always a solution; hence, we assume that $y\neq(1-b)/(az)$ for the rest of the proof. We have the recurrence relations as follows: $\theta_{-1}=(1-b)/a,~\theta_0=0$ and
\begin{subequations}
\label{eq:fde:2}
{\begin{align}
\label{eq:fde:3}
(2k+2(1-b))\theta_{2k}+a\sum_{j=0}^{2k-1}\theta_{2k-1-j}\theta_j=0,\\
\label{eq:fde:4}
(2k+1+2(1-b))\theta_{2k+1}+a\sum_{j=0}^{2k}\theta_{2k-j}\theta_j+b\theta_k=0.
\end{align}}
\end{subequations}
If $n+2(1-b)\neq0,\forall n\geq1$, then $\theta_n=0,~n\geq0$ can be found recursively from (\ref{eq:fde:2}) and we again obtain $y=(1-b)/(az)$. Hence, we must have $N+2(1-b)=0$ for some $N\in\mathbb{N}$. In this case, we have $\theta_{0}=\cdots=\theta_{N-1}=0$ and $\theta_N$ is a parameter. When $\theta_N=0$, we obtain $\theta_n=0$ for all $n\geq0$, and $y=(1-b)/(az)$. Thus, we assume that $\theta_N\neq0$. Similar arguments show that $y$ admits only a pole at $z=0$. If $y$ is rational, then we obtain $\theta_k=0$ for all large $k$ and $y=P(z)/z$ for some polynomial $P$. Then
(\ref{eq:intro:7}) reduces to
\begin{equation}
\label{eq:fde:5}
zP'-P+aP^2+bP(z^2)=0.
\end{equation}
Suppose that $P(z)$ is a polynomial of degree $p$ and so $P(z)=\sum_{j=-1}^{p-1}\theta_{j}z^{j+1}$, which by substituting into (\ref{eq:fde:5}) yields $\theta_{p-1}=-b/a$ and $\theta_{-1}=(1-b)/a$. We also find from the coefficients of $z^j$ in (\ref{eq:fde:5}) for $j=2p,2p-1,\dots,p+1$ subsequently that $\theta_{i-1}=0$ for $i=p-1,p-2,\dots,1$. Moreover, the coefficient of $z^p$ in (\ref{eq:fde:5}) shows that $2(1-b)+1-p=0$ implies $N=p-1$. Hence, $y=(1-b)/(az)-bz^{N}/a$ and this corresponds to $\theta_N=\theta_{p-1}=-b/a$. For $\theta_N\neq 0,~-b/a$, we have $y$ as a transcendental function. Theorems \ref{thm:nev:9} and \ref{thm:nev:10} (and the last part of the proof of the latter) assert that $y$ cannot be meromorphic and must admit $C(0,1)$ as a natural boundary either as a single-valued function defining on $D(0,1)$ or a Riemann surface lying over $D(0,1)$ with the exception of some discrete points on each Riemann sheet such that $C(0,1)$ is a natural boundary of each Riemann sheet.

Now we prove the last statement. If $\tilde{r}=1$, we are done. Assume that $\tilde{r}<1$ and there exists a single-valued meromorphic $y$ defined on $D(0,1)$. In this case, $y$ has a pole at some point $z_0\in C(0,\tilde{r})$, and the pole must be simple; otherwise, if $z_0$ is a pole of order $\ell>1$, then $z_0^2$ is a pole of order $2\ell$ contradicting the assumption of $\tilde{r}$. Similarly, if there is a pole at $z_0'$ of $y$ outside $\overline{D}(0,\tilde{r})$, the pole must be simple; otherwise, we obtain a sequence of poles $\{(z_0')^{2^n}\}$ accumulating to $0$, a contradiction. Hence, all poles of $y$ are simple. Consider the expansion of $y$ at $z_0$, which is $y(z)=c_{-1}/(z-z_0)+O(1)$. Then $y(z^2)=c_{-1}/[2\sqrt{z_0}(z-\sqrt{z_0})]+O(1)$.
Assume that the expansion of $y$ at $\sqrt{z_0}$ is $y(z)=d_{-1}/(z-\sqrt{z_0})+d_0+O(z-\sqrt{z_0})$. From (\ref{eq:intro:7}), we obtain $d_{-1}=1/a$ and $d_0=-bc_{-1}/(4\sqrt{z_0})$, which shows that $y'+ay^2=-bc_{-1}/[2\sqrt{z_0}(z-\sqrt{z_0})]+O(1)$. Let $N^*(r,\cdot)$ denote the counting function of the poles in $0<|z|\leq r$, so that the pole at the origin is excluded; a single point does not affect the argument. Then $N^*(r,y'+ay^2)=N^*(r,y)$, and since $N^*(r,-by(z^2))=N^*(r^2,y)$, we obtain $N^*(r^2,y)=N^*(r,y)$ for $0\leq r<1$. Hence $N^*(r,y)=N^*(r^{2^n},y)\to0$ as $n\to\infty$ for every $0<r<1$, so $y$ has no poles in $0<|z|<1$; that is, $y$ only has a pole at $z=0$. Therefore $y$ must have a branch point on $C(0,\tilde{r})$.
\end{proof}

Theorem \ref{thm:fde:1} shows that $y$ has $C(0,1)$ as a natural boundary. Statement (5) shows that a transcendental function $y$ of the form (\ref{eq:fde:1}) is either single-valued and meromorphic on $D(0,1)$, or multivalued, according to whether $\tilde{r}=1$ or $\tilde{r}<1$. When $\tilde{r}<1$, we do not know whether the number of Riemann sheets is finite or infinite, and we leave this problem for future work. We now prove that both cases $\tilde{r}=1$ and $\tilde{r}<1$ occur for each $a\neq0$ and $b=(N+2)/2$, $N\in\mathbb{N}$.

\subsection{Existence of the case $\tilde{r}=1$}
\label{subsec:fde:1}

We recall here that $b=(N+2)/2$ and $\theta_N$ is a parameter different from $0$ and $-b/a$. Setting $\beta_j=-a\theta_j$ in (\ref{eq:fde:2}), we obtain
\begin{subequations}
\label{eq:fde:6}
\begin{align}
\label{eq:fde:7}
&\beta_{2k}=\frac{1}{2k-N}\sum_{j=0}^{2k-1}\beta_{2k-1-j}\beta_j,\\
\label{eq:fde:8}
&\beta_{2k+1}=\frac{1}{2k+1-N}\sum_{j=0}^{2k}\beta_{2k-j}\beta_j-\frac{N+2}{2(2k+1-N)}\beta_k.
\end{align}
\end{subequations}
Since $\theta_i=0$ for each $i=0,\dots,N-1$, then the term $\beta_j$ that could be non-zero is those of the indices $j=i(N+1)-1$, for $i\geq0$.
We take the parameter $\theta_N$ so that $|\beta_N|\leq1$. For $N\geq2$, we will prove by induction that $|\beta_{i(N+1)-1}|\leq|\beta_N|$ for $i\geq1$. By the hypothesis of induction, we obtain from (\ref{eq:fde:6}) the following
\begin{subequations}
\begin{align*}
|\beta_{i(N+1)-1}|&\leq\frac{1}{i(N+1)}\sum_{j=1}^i|\beta_{j(N+1)-1}\beta_{(i-j+1)(N+1)-1}|+\frac{N+2}{2(iN+i)}|\beta_N|\\
&\leq\frac{i|\beta_N|^2}{i(N+1)}+\frac{N+2}{2(N+1)i}|\beta_N|\leq\frac{N+4}{2(N+1)}|\beta_N|\leq|\beta_N|.
\end{align*}
\end{subequations}
This shows that $\sum_{i=0}^\infty |\beta_i||z|^i\leq|\beta_N|\sum_{i=0}^\infty |z|^i$, which converges for all $|z|<1$. For $N=1$, we alter the above argument slightly. We find that $|\beta_3|\leq|\beta_1|^2/2+3|\beta_1|/4\leq5|\beta_1|/4$ and $|\beta_5|=|\beta_1g_3|/2\leq5|\beta_1|/8$. Similarly, we have $|\beta_7|\leq|\beta_1|$. For $i\geq4$, we will prove by induction that $|\beta_{2i-1}|\leq|\beta_1|$. From (\ref{eq:fde:6}), we obtain
\begin{align*}
|\beta_{2i+1}|&\leq\frac{1}{2i}\sum_{j=1}^i|\beta_{2j-1}||\beta_{2(i-j+1)-1}|+\frac{3}{4i}|\beta_1|\\
&\leq\frac{(i-4)|\beta_1|^2+5|\beta_1|^2/2+5|\beta_1|^2/4}{2i}+\frac{3}{4i}|\beta_1|\leq\frac{4i+5}{8i}|\beta_1|\leq|\beta_1|.
\end{align*}
This again shows that $\tilde{r}\geq1$, and hence $\tilde{r}=1$ by Theorem \ref{thm:fde:1}(4).

\subsection{Existence of the case $\tilde{r}<1$}
\label{subsec:fde:2}

We will choose $\theta_N$ such that $\beta_j$ is positive and satisfies $\beta_N\geq\alpha\xi N+4(N+2)$ , where $\xi>1$ is a fixed positive number and $\alpha\geq7/2$. We will prove by induction that $\beta_{i(N+1)-1}\geq\alpha\xi^iiN+4(N+2)$. First, we find that
\begin{subequations}
\begin{align*}
\beta_{2N+1}&=\frac{1}{N+1}\beta_N^2-\frac{N+2}{2(N+1)}\beta_N=\frac{\beta_N}{N+1}\left(\beta_N-\frac{N+2}{2}\right)\\
&\geq\frac{1}{N+1}(\alpha\xi N+3(N+2))^2\geq\frac{\alpha^2\xi^2N^2}{N+1}+\frac{3\alpha(N+2)N}{N+1}\xi\geq2\alpha\xi^2N+4(N+2).
\end{align*}
\end{subequations}
Let $G(i)=\alpha\xi^iiN$. For $i\geq2$, we have
\begin{equation*}
\label{eq:fde:9}
\begin{aligned}
\beta_{(i+1)(N+1)-1}\geq\frac{1}{iN+i}\sum_{j=1}^i\beta_{j(N+1)-1}\beta_{(i-j+1)(N+1)-1}-\frac{\delta(N+2)}{i(N+1)}\beta_{\frac{(i+1)(N+1)-1}{2}},
\end{aligned}
\end{equation*}
where $\delta=1$ if $[(i+1)(N+1)-1]/2$ is a positive integer, otherwise $0$. Since $\beta_{i}^2-(N+2)\beta_i\geq(G(i)+3(N+2))^2$, then
\begin{subequations}
\begin{align*}
\beta_{(i+1)(N+1)-1}&\geq\frac{1}{iN+i}\sum_{j=1}^i(G(j)+3(N+2))(G(i+1-j)+3(N+2))\\
&\geq\frac{\alpha^2\xi^{i+1}N^2}{iN+i}\sum_{j=1}^ij(i+1-j)+\frac{3\alpha(N+2)\xi iN}{iN+i}\\
&\geq\frac{\alpha^2\xi^{i+1}N^2i(i+1)(i+2)}{6i(N+1)}+\frac{21(N+2)}{4}\geq(i+1)\alpha N\xi^{i+1}+4(N+2),
\end{align*}
\end{subequations}
where the last inequality follows from $\alpha N(i+2)\geq4\alpha N\geq 6N+6$. This shows that $\sum_{i=1}^\infty \beta_jz^{j+1}=\sum_{i=1}^\infty \beta_{i(N+1)-1}(z^{N+1})^i\geq \alpha\xi N\sum_{i=1}^\infty (i+1)(\xi z^{N+1})^i$. The series $\sum_{i=1}^\infty (i+1)(\xi z^{N+1})^i$ has the radius of convergence $\xi^{-1/(N+1)}<1$. Hence, $\tilde{r}<1$.

\begin{remark}
\label{rem:fde:1}
When $y$ is analytic at $0$ (the solution given by \cite[Th.2.1]{oberg1971local}), then we set $y=\sum_{i=0}^\infty\theta_iz^i$, where $\theta_0$ is a parameter. It can be easily checked that the only rational solutions here are $0$ and $-b/a$. Thus, all other solutions are transcendental. Moreover, there exists a transcendental solution even if $b\neq (N+2)/2$ for any $N\in\mathbb{N}$. The conclusions for the transcendental solutions here are similar to Theorem \ref{thm:fde:1}, those are: when $a=0$, the transcendental solutions are single-valued defined on $D(0,1)$. For $a\neq0$, the transcendental solutions are either single-valued or multi-valued. All solutions possess $C(0,1)$ as a natural boundary. The proof is more or less similar to the proof of Theorem \ref{thm:fde:1}. For the discussion of $\tilde{r}=1$ and $\tilde{r}<1$ when $a\neq0$, the proof for those $b=(N+2)/2$ and, more generally, for $b>0$ are similar to the above discussions to show both the existence of the case $\tilde{r}=1$ and $\tilde{r}<1$. For the remaining values of $b$, we have not yet obtained this result.
\end{remark}

\section{Singularity structure of special solutions of the q-DE (\ref{eq:intro:8})}
\label{sec:qde}

We now discuss the special solutions of (\ref{eq:intro:8}) that have a pole at $z=0$, where $q\neq0$ is fixed. Theorem \ref{thm:nev:8} and Corollary \ref{cor:nev:1} show that such solutions exist when $|q|\leq1$. The analytic continuation is trivial when $|q|<1$, but requires some discussion when $|q|=1$. The existence of such a solution when $|q|>1$ is proved below. The most interesting case is $|q|>1$, where we show that a transcendental solution must admit branch points. Recall that we seek a solution $y$ with a pole at $0$. Balancing the dominant terms of a Laurent ansatz shows that the pole must be simple: if $y\sim cz^{-p}$ with $p\geq2$, then $by^2(qz)$ is the only term of order $z^{-2p}$, which forces $c=0$. Hence $y$ satisfies
\begin{equation}
\label{eq:qde:1}
y=\frac{\theta_{-1}}{z}+\sum_{i=0}^\infty\theta_iz^i,~\theta_{-1}\neq0.
\end{equation}
Our main results are the following theorems.

\begin{theorem}
\label{thm:qde:1}
Let (\ref{eq:qde:1}) be a formal series solution of (\ref{eq:intro:8}). If $n+1+2q^{n+2}\neq0$ for all $n\geq0$, then (\ref{eq:qde:1}) converges in a punctured neighbourhood of $0$. If $n+1+2q^{n+2}=0$ for some $n=N$, then $N$ is the only non-negative integer root, and, under the resonance condition, the series (\ref{eq:qde:1}) again converges in a punctured neighbourhood of $0$.
\end{theorem}

\begin{proof}
Since the case $|q|\leq1$ is covered by Theorem \ref{thm:nev:8} and Corollary \ref{cor:nev:1}, it suffices to consider $|q|>1$. Substituting (\ref{eq:qde:1}) into (\ref{eq:intro:8}), we obtain
\begin{equation}
\label{eq:qde:2}
\theta_{-1}=\frac{q^2}{b},~\theta_0=\frac{-aq}{2b},~(n+1+2q^{n+2})\theta_{n+1}+a\theta_n+b\sum_{j=0}^n\theta_j\theta_{n-j}q^n=0,
\end{equation}
for $n\geq0$. Setting $\beta_j=b\theta_j$, we obtain
\begin{equation}
\label{eq:qde:3}
(n+1+2q^{n+2})\beta_{n+1}+a\beta_n+\sum_{j=0}^n\beta_j\beta_{n-j}q^n=0.
\end{equation}
First, we assume that $n+1+2q^{n+2}\neq0$ for all $n\geq0$. Let $V_q=\sup_{n\geq0}\{1,|2+(n+1)/q^{n+2}|^{-1}\}$, which is finite. We prove by induction that $|\beta_j|\leq|\beta_0|C_{j}|a|^jV_q^j|q|^{-j}$. Here, $C_{n}$, $n\geq0$, is the $n$-th Catalan number, defined by $C_{0}=1$ and $C_{n}:=\binom{2n+1}{n}/(2n+1)$; see \cite[Ch.5]{koshy2008catalan}. It is also known that
$\sum_{j=0}^nC_{j}C_{n-j}=C_{n+1}$, which shows that $C_{n}$ is an increasing sequence. For $n=0$, we obtain
\begin{equation*}
\label{eq:qde:4}
\begin{aligned}
|\beta_1|=\frac{1}{2|q|^2}|\beta_0|\left(|a|+\frac{|a||q|}{2}\right)\leq\frac{3|\beta_0|}{4|q|}V_q|a|\leq\frac{|a||\beta_0|}{|q|}V_qC(1).
\end{aligned}
\end{equation*}
From (\ref{eq:qde:3}), we have for $n\geq1$
\begin{align*}
|\beta_{n+1}|&\leq\frac{1}{|q|^{n+2}}V_q\left(|a|C_{n}V_q^n|a|^n|\beta_0||q|^{-n}+|a^n|V_q^nC(n+1)|\beta_0|^2 \right)\\
&\leq\frac{V_q^{n+1}|a|^{n+1}|\beta_0|}{|q|^{n+1}}\left(\frac{C_{n}}{|q|}+\frac{C_{n+1}}{2} \right)\\
&\leq\frac{V_q^{n+1}|a|^{n+1}|\beta_0|}{|q|^{n+1}}C_{n+1},
\end{align*}
where we use $C_{n+1}=2(2n+1)C_{n}/(n+2)\geq2C_{n}$ for $n\geq1$.
Hence, $\sum_{j=0}^\infty|\beta_jz^j|\leq|\beta_0|\sum_{j=0}^\infty |aV_qz/q|^jC_j$ converges within $|4aV_qz/q|<1$ \cite[p.120]{koshy2008catalan}.

When $n+1+2q^{n+2}=0$, we will show at the end that there exists at most one non-negative integer root $N$. Under the assumption of the resonance condition, if there is still a formal series (\ref{eq:qde:1}), we will show that it is still a Laurent series. Here, $\theta_N$ is a parameter. We only need to alter the definition of $V_q$ to $V_q=\sup_{n\geq0,~n\neq N}\{1,|2+(n+1)/q^{n+2}|^{-1},|\beta_N||q|/(|\beta_0a|)\}$, which shows that $|\beta_N|\leq |\beta_N|C_{N}\leq|\beta_0||aV_q|^NC_N|q|^{-N}$ and the proof above remains valid.

Now we show that $N$ is unique. First we have $0$ is a root iff $2q^2+1=0$, and $1$ is a root iff $q^3+1=0$. For $|q|>1$, we may assume that $n\geq2$. Suppose that $N$ and $N+N_1$ are both positive integer roots. We find that $N+1+2q^{N+2}=N+N_1+1+2q^{N+N_1+2}=0$, which leads to $(N+N_1+1)^{N+2}2^{N_1}+(-(N+1))^{N+N_1+2}=0$.
So both $N+1$ and $N_1$ must be even. Let $N+1=2p$ and $N_1=2k$, then we have $(p+k)^{2p+1}=p^{2p+2k+1}$ which is equivalent to $F(p)=F(p+k)$ for a function $F(x):=\ln x/(2x+1),~x\geq2$. However, $F'(x)<0$ for $x\geq3$. We checked numerically that $F(n_1)\neq F(n_2)$ for distinct positive integers $n_1,n_2\in[2,6]$ and that $F(2)>F(6)$. Hence, if a positive integer root $N$ exists, it is unique.
\end{proof}

\begin{theorem}
\label{thm:qde:2}

Let $y$ be the maximal analytic continuation of the Laurent series solution (\ref{eq:qde:1}) of (\ref{eq:intro:8}). Then we have the following.
\begin{enumerate}
    \item If $|q|\neq1$ and $a\neq0$, then a non-constant rational solution occurs only for $q=2$, namely $y=(4-az)/(bz)$; moreover, for $q=2$ the solution is always $y=(4-az)/(bz)$.
    \item If $2q^2+1=0$, then $a=0$, and $y$ is either $\theta_{-1}/z$ or a transcendental meromorphic function. For all other $q$ with $|q|<1$, and $a\neq0$, $y$ is transcendental and meromorphic.  These solutions have only one pole, at $z=0$.
    \item Suppose that $a=0$. If $n+1+2q^{n+2}\neq0$ for all $n\geq2$, then $y=\theta_{-1}/z=q^2/(bz)$ is rational. If $n+1+2q^{n+2}=0$ for some $n=N\neq1$, then $y$ is either $\theta_{-1}/z$ or transcendental (and meromorphic if $|q|<1$). The case $N=1$ is included in case (4).
    \item Let $a$ be arbitrary, and suppose that $q$ is a primitive $m$-th root of unity. Then (\ref{eq:intro:8}) reduces to a quadratic system in $m$ variables. If $m=2$, the general solutions are in the class $W$, and so is $y$. If $m\geq3$ is odd, every meromorphic solution is in the class $W$. If $m=6$, then $a=0$ and $y=\theta_{-1}/z$. The class $W$ is defined, and this case is proved, in Section \ref{subsec:qde:2}.
    \item When $|q|=1$, $a\neq0$ and $q$ is not a root of unity, then $y$ is meromorphic on $D(0,r_q)$ with $C(0,r_q)$ as a natural boundary, where $r_q$ is the radius of convergence of the analytic part of (\ref{eq:qde:1}).
    \item If $|q|>1$, $q\neq2$, $a\neq0$ and $n+1+2q^{n+2}\neq0$ for all $n$, then $y$ is transcendental but not meromorphic.
\end{enumerate}
\end{theorem}

In case 4 above, when $m\geq3$ and $m\neq6$, we do not know whether $y$ is transcendental or meromorphic.
In case 6, $y$ defines a Riemann surface over $\mathbb{C}$, whose branch points are described in Section \ref{subsec:qde:1}. When $n+1+2q^{n+2}=0$ for some $n$, the situation is the same; however, we conjecture that $a=0$ there.

\begin{proof}
Now we prove Case 1. We recall $|q|\neq1$, then it can be checked that if $z_0\neq0$ is a pole of $y$, then we obtain a sequence of poles $\{z_0q^n,n\geq0\}$ of $y$. Thus, if $y$ is rational, then there is only one pole at $z=0$. Let $y=P(z)/z$ for some polynomial $P$, then by substituting into (\ref{eq:intro:8}), we obtain $zP'-P+azP+bP^2(qz)/q^2=0$. If $\deg P=p$, then it is easy to see that $p=1$. Thus, $y=\theta_{-1}/z+\theta_0$. From the recurrence relation (\ref{eq:qde:2}) shows that $q(2q^2+1)\theta_1=b\theta_0^2(2-q)$. If $2q^2+1=0$, then we obtain $\theta_0=0$ and so $a=0$, i.e., $y=q^2/(bz)$. If $2q^2+1\neq0$, then $\theta_0=0$ and so $a=0$ where we obtain the same rational solution. We could obtain a different rational solution which can be seen by setting $q=2$, and we find that $y=\theta_{-1}/z+\theta_0=(4-az)/(bz)$. It is easy to see that for $q=2$, we obtain $\theta_i=0$ for $i\geq1$ and so $y=(4-az)/(bz)$.

In Case 2, when $2q^2+1=0$ yields a resonance condition $a=0$ and $\theta_1$ is free. If $\theta_1=0$, we find that $\theta_i=0$ and $y=\theta_{-1}/z$. For $\theta_1\neq0$, Case 1 shows that $y$ is transcendental. For the rest of $|q|<1$, we obtain $n+1+2q^{n+2}\neq0$ for all $n$. Case 1 again shows that $y$ is transcendental. Since $|q|<1$, $y$ is meromorphic. It is easy to see that $z=0$ is the only pole of $y$ by arguing similarly as in Case 1.

When $a=0$, we can find subsequently from (\ref{eq:qde:2}) $\theta_i=0$ for $i\geq0$ unless $n+1+2q^{n+2}=0$ for some $n$. If $N\neq1$ is a root of $n+1+2q^{n+2}=0$, then $\theta_N$ is free. If $\theta_N=0$, we again obtain $y=\theta_{-1}/z$; otherwise, $y$ is transcendental following Case $1$. When $N=1$, we obtain $q=-\omega$, for $\omega$ any root satisfying $\omega^3=1$, and we include in Case 4.

Case 4 is proved in Section \ref{subsec:qde:2}.

When $|q|=1$ and not a root of unity, then arguing similar to Case 1 shows that any meromorphic solution $y$ admits only one pole. In this case, we also obtain $2m(r,y)=m(r,y^2(qz))\leq m(r,y)+S(r,y)$. Hence, the only meromorphic solution here is the rational solution with one pole. We again will only obtain $y=\theta_{-1}/z$ and $a=0$. For $a\neq0$, the solution is not meromorphic. This means that the analytic part of (\ref{eq:qde:1}) admits a radius of convergence $r_q$ with at least a singular point $z_1\in C(0,r_q)$. However, the dense set of singular points $\{z_1q^n,~n\geq0\}$ shows that $C(0,r_q)$ is a natural boundary. (If $y$ is singular at $z_1$, then it is singular at $qz_1$: otherwise $y^2(qz)$ would be analytic near $z_1$, and then $y$ would be analytic near $z_1$ as well, because $y'+ay=-by^2(qz)$ is a linear first-order equation for $y$.)

Finally, we arrive at the case $|q|>1$. Based on Case 2, we assume that $q\neq2$ and $a\neq0$. We show that (\ref{eq:qde:1}) is not meromorphic. Here, $y$ has either one pole (at $z=0$) or infinitely many poles as shown in Case 1. In the latter case $y$ is not meromorphic. This is because if $z_1\neq0$ is a pole, then the sequence of poles $q^{-n}z_1$ converging to $0$. Hence, $y$ has only one pole. Moreover, $2m(r,y)\leq2m(r|q|,y)=m(r,y^2(qz))\leq m(r,y)+S(r,y)$; thus, $m(r,y)=S(r,y)$. The case where $y$ has only one pole is the rational solution mentioned in Case 1 corresponding to $a=0$. Hence, when $a\neq0$, $y$ is not meromorphic, and it must admit a branch point on $C(0,r_q)$, where $r_q$ is the radius of convergence of the analytic part of (\ref{eq:qde:1}). More precisely, let $w_0\in C(0,r_q)$ be a singular point of $y$. Since $|q|>1$, the point $w_0/q$ lies inside the disc of convergence, so $h(z):=-(y'(z)+ay(z))/b$ is analytic near $w_0/q$, and (\ref{eq:intro:8}) gives $y^2(qz)=h(z)$ there. Hence $y(\zeta)=\pm\sqrt{h(\zeta/q)}$ near $\zeta=w_0$, which is analytic unless $h$ has a zero of odd order at $w_0/q$. Therefore every singular point of $y$ on $C(0,r_q)$ is a quadratic branch point at which $y'+ay$ has a zero of odd order at $w_0/q$.  We end this case with a remark on the resonance condition of the resonance $N$ of $N+1+2q^{N+2}=0$. It is easy to find that $\theta_i=b^i\theta_0^{i+1}R_i(q),~i<N$, where $R_i$ are rational functions. The resonance condition obtaining from (\ref{eq:qde:2}) admits the form $b^{N+1}\theta_{0}^{N+2}P_N(q)=0$ for some polynomial $P_N$. If $P_N=0$, we would obtain $y$ as a transcendental solution. However, we believe that $P_N\neq0$ for all $N$ and this leads to $\theta_0=a=0$; the numerical result in Mathematica for $N\leq15$ shows that it is true.
\end{proof}

\subsection{The nature of the branch points}
\label{subsec:qde:1}
We now describe the structure of the branch points of the transcendental solution $y$ when $|q|>1$ and $q\neq2$. As shown in Theorem \ref{thm:qde:2}, $y$ must admit a branch point on $C(0,r_q)$, the circle of convergence of the analytic part of (\ref{eq:qde:1}). Let $w_0\in C(0,r_q)$ be a branch point (there may be more than one on the circle). Then $y'+ay$ has a zero of odd order at $w_0/q$, and $w_0$ is a quadratic branch point of $y$ (see the proof of Theorem \ref{thm:qde:2}, case 6). We denote $y_0$ to be the sheet containing (\ref{eq:qde:1}) and $y_1$ to be the other Riemann sheet connecting with $y_0$ at $w_0$. Moreover, $L[y_0]=y_i^2(qz)$ where $L[y]:=-(y'(z)+ay(z))/b$. Then on $y_i$, we have a branch point of fourth order at $qw_0$; we do not know whether $y_0$ and $y_1$ are glued through the branch point $qw_0$. To see the fourth-order branching at $qw_0$, note that in the generic case where $y'+ay$ has a simple zero at $w_0/q$, we have $y_0(z)=y_0(w_0)+c(z-w_0)^{1/2}+\cdots$ with $c\neq0$, so $y_0'(z)\sim\frac{c}{2}(z-w_0)^{-1/2}$. Hence $y_i^2(qz)=L[y_0](z)\sim-\frac{c}{2b}(z-w_0)^{-1/2}$, and so $y_i(\zeta)\sim C(\zeta-qw_0)^{-1/4}$ near $\zeta=qw_0$ for some constant $C\neq0$; in particular $y_i\to\infty$ there. Repeating the arguments, we obtain that at $q^nw_0$, $y_i$ is branched of order $2^{n+1}$. Of course, there may be other branch points that are not inherited from $w_0$. There could be a branch point at $w_1$ (a new stem) on some sheet, as the first appearance of a new type of branch point, similar to the examples discussed in Section \ref{subsec:poinc:2}; here it corresponds to an odd-order zero different from the zero corresponding to $w_0$.

\subsection{Meromorphic solutions when $q$ is a root of unity}
\label{subsec:qde:2}
In this subsection, we prove case (4) of Theorem \ref{thm:qde:2}. Class $W$ consists of the elliptic functions and their degenerations, that is, rational functions of $e^{az}$ ($a\in\mathbb{C}$) and rational functions of $z$ (see \cite{eremenko1982meromorphic}). A meromorphic solution of an algebraic differential equation has the \emph{finiteness property} if there are only finitely many distinct formal Laurent series solutions at its poles (up to translation of the pole). 

Let $q$ be a primitive $m$-th root of unity with $m\neq6$. Then (\ref{eq:intro:8}) can be written as a quadratic system in the variables $x_i(z)=y(q^iz)$, $i=0,\dots,m-1$ (with indices modulo $m$):
\begin{equation}
\label{eq:qde:5}
    x_i'(z)/q^i+ax_i(z)+bx_{i+1}^2(z)=0.
\end{equation}
If $y$ is meromorphic, then by definition $T(r,x_i)=T(r,x_j)$ and $m(r,x_i)=m(r,x_j)$ for all $i,~j$. Moreover, (\ref{eq:intro:8}) shows that $m(r,y)=S(r,y)$ and so $m(r,x_i)=S(r,x_i)$ for all $i$. If $y$ has only one pole at $z=0$, then we obtain $T(r,y)=\ln r+S(r,y) $. This means that $y=\theta_{-1}/z+\theta_0$ and by substituting into (\ref{eq:intro:8}), we find that it occurs if and only if $a=0$ with $y=q^2/(bz)$. When $a\neq0$, we cannot say more about the transcendence or meromorphicity of $y$ in general. However, the Laurent expansions of the solutions of (\ref{eq:qde:5}) at poles $z_0\neq0$ give the following information on the resonances. Writing $x_i=x_{0,i}/t+x_{1,i}+\cdots$ with $t=z-z_0$, we obtain
\begin{equation}
\label{eq:qde:6}
    x_{0,i}=q^{i}bx_{0,i+1}^2,\qquad \prod_{i=0}^{m-1}x_{0,i}=q^{-\frac{m(m-1)}{2}}b^{-m}.
\end{equation}
Using (\ref{eq:qde:6}), we find that the coefficient matrix of the $x_{n,i}$ has determinant $q^{-\frac{m(m-1)}{2}}\big((n-1)^m+(-1)^{m-1}2^m\big)$. When $m$ is odd, then there are no non-negative integer resonances. Hence, the expansions around the poles of $x_i$ are unique; so $x_i$ has the finiteness property, and \cite[Theorem 3, third case]{eremenko1982meromorphic} shows that the meromorphic solutions\footnote{This does not imply that the solution $y$ is meromorphic.} $x_i$ are in the class $W$; see also \cite[Lemma 3.5]{meng2024meromorphicsolutions3dlotkavolterra}. However, when $m$ is even, there is exactly one positive integer resonance at $n=3$ and we believe that we can construct a first integral
of (\ref{eq:qde:5}). We have so far constructed only the case $m=2$ as for the higher dimensional cases require more work and we leave this for future work. Now we investigate the case $m=2$ (so $q=-1$). The conserved quantity is given by: $x_0^3+x_1^3+3ax_0x_1/b=g$. The Painlev\'e test shows that there are three types of pole, each of which is simple. The residues of the poles for $x_0$ and $x_1$ are $\omega b^{-1},$ and $-\omega^2 b^{-1}$ respectively, where $\omega$ can be any roots of $\omega^3=1$. The general solution of $x_0$ and $x_1$ can be found from the ans\"atze in \cite[§7.5.2]{armitage2006elliptic}. We set a parameter corresponding to the movable poles so that the solution $y$ coincides with (\ref{eq:qde:1}). We found the following cases:
\begin{equation*}
(y,a)=\left(\frac{1}{bz},0\right),~\left(\frac{\pm\mathrm{i}\sqrt{3}a}{b(e^{\pm\mathrm{i}\sqrt{3}a z}-1)}+\frac{a}{2b}\pm\frac{\mathrm{i}\sqrt{3}a}{2b},a\right),
\end{equation*}
corresponding to $g=a^3/b^3$, for $g\neq a^3/b^3$, $y$ is an elliptic function (we omit the expression).
Finally, when $m=6$, then $q=-\omega$ and $n+1+2q^{n+2}=0$ at $n=1$ with the resonance condition $a\theta_0=0$. Thus, we obtain $a=0$ and $y=\theta_{-1}/z$.

\section{Discussion}
\label{sec:disc}

In Section \ref{subsec:poinc:2}, we saw various examples of the singularity structure of the analytic continuation of the Poincar\'e function, equivalently of the solutions of (\ref{eq:intro:3}). Notably, a stem can first occur on a sheet of the $k$-th generation for some polynomials of degree greater than $3$. So far, we have not been able to construct such an example of degree $3$, and it would be interesting to know whether this is impossible. Our examples of this delayed appearance are also rather simple. We would like to construct an example combining Example \ref{ex:poinc:2} with Example \ref{ex:poinc:6} or \ref{ex:poinc:7}, in which critical points are mapped to other critical points and the appearance of the stem is delayed. We plan to address this in future work.

For equation (\ref{eq:intro:7}), our results leave a gap for certain values of the parameter $\theta_N$, equivalently of $\beta_N=-a\theta_N$, namely those not covered by Sections \ref{subsec:fde:1} ($|\beta_N|\leq1$) and \ref{subsec:fde:2} ($\beta_N\geq\alpha\xi N+4(N+2)$), and for some $a$ and $b$.  We believe that this gap can be closed. When the radius of convergence is less than $1$, we would like to see if the Riemann surface has infinite sheets; we expect it to be infinite. Another interesting problem in that case is to determine all the stems lying on the sheet $y_0$, and, even more interestingly, whether a delayed stem similar to those in Examples \ref{ex:poinc:6} and \ref{ex:poinc:7} can occur. The same questions can be asked for equation (\ref{eq:intro:8}) when $|q|>1$.

In Section \ref{sec:fde}, we only considered solutions of (\ref{eq:intro:7}) with a simple pole at $0$. Unlike Mahler's equation (\ref{eq:intro:5}), where poles of higher order can be reduced to the simple case by composing with $z^p$ \cite{mahler1981special}, a Laurent ansatz $y=cz^{-p}+\cdots$ with $p\geq2$ for (\ref{eq:intro:7}) leads to a recurrence relation whose formal solution has factorially growing coefficients, so the convergence of the analytic part is unclear. Once such convergence is established, we expect the singularity structure to be described as in the case of a simple pole, similarly to Mahler's equation \cite{mahler1976class,mahler1981special}.

Several questions about (\ref{eq:intro:8}) also remain open. When $|q|>1$ and $N+1+2q^{N+2}=0$ for some $N$, the resonance condition takes the form $b^{N+1}\theta_0^{N+2}P_N(q)=0$ for a polynomial $P_N$ (see the proof of Theorem \ref{thm:qde:2}). We conjecture that $P_N\neq0$ for all $N$, so that the resonance forces $a=0$; we have verified this numerically for $N\leq15$. When $q$ is a root of unity of even order $m\geq4$, we expect the system (\ref{eq:qde:5}) to have a first integral associated with the resonance at $n=3$, as in the case $m=2$, but we have not been able to construct it. Finally, we have not described how the Riemann sheets are glued, either for the Poincar\'e functions in Section \ref{sec:poinc} or for the solutions in Sections \ref{sec:fde} and \ref{sec:qde}; for the quadratic family (\ref{eq:intro:4}), this was done in \cite[§4]{halburd2025extensionanalyticsolutionsfirstorder}.

\section*{Declaration of AI Use}
Most of the results presented in this paper were obtained, without the use of AI, as part of the author's PhD thesis \cite[Ch.5,6]{meng2026application} at University College London. DeepSeek was used to find the example (\ref{eq:poinc:18}), expansion of Puiseux series, and to help with the asymptotics of $S_r$ in Example \ref{ex:poinc:6}. Claude (Anthropic) was used to improve the language and presentation, and to check several errors of Examples \ref{ex:poinc:1}--\ref{ex:poinc:5}. It also pointed out that the proof of Theorem \ref{thm:nev:9} in \cite{meng2026application} omitted the case $\rho=0$, for which the author had obtained a proof in an unpublished note but had inadvertently left it out of the thesis. It also highlighted an issue concerning the condition $o(T(r,y))$.

\section*{Acknowledgements}
The author would like to thank Professor Rod Halburd and Dr.~Bishal Deb for useful comments, and Professor Rod Halburd and Professor Qingping Liu for their support and encouragement. This work was supported by the National Natural Science Foundation of China (NNSFC) (Grant Nos.~12571272 and 12171474).

\appendix
\section{Proofs of Theorem \ref{thm:nev:8} and Corollary \ref{cor:nev:1}}
\label{sec:appA}

Here we prove Theorem \ref{thm:nev:8} and Corollary \ref{cor:nev:1}. 

\begin{proof}[Proof of Theorem \ref{thm:nev:8}]
We begin with (\ref{eq:nev:3}), and first assume that $h\equiv1$. Then the recurrence relation (\ref{eq:nev:5}) becomes
\begin{equation}
\label{eq:appA:1}
(n+1-2\theta_{-1}a_0)\theta_{n+1}=P_n(\theta_{-1},\theta_0,\dots,\theta_{n}),
\end{equation}
where $P_{-1}=(a_1\theta_{-1}+b_0+d_1\theta_{-1})\theta_{-1}$ and the $P_n$ are polynomials, whose explicit form we do not need. Since $n+1-2\theta_{-1}a_0\neq0$ for all $n$, there exists $K>0$ such that $|n+1-2\theta_{-1}a_0|^{-1}\leq K$ for $n\geq0$. Let $\hat{y}=\sum_{j=0}^\infty\hat{\theta}_jz^j$ be defined by $\hat{\theta}_0=|\theta_0|$ and
\begin{equation}
\label{eq:appA:2}
\hat{y}=Kz(\hat{a}\hat{y}^2+\hat{b}\hat{y}+\hat{d}\hat{y}+\hat{f})+\hat{\theta}_0,
\end{equation}
where
\begin{align*}
&\hat{a}=\sum_{j=0}^\infty|a_j|z^j,\qquad \hat{b}=\sum_{j=0}^\infty\big(|b_j|+|2\theta_{-1}a_{j+1}|\big)z^j,\qquad \hat{d}=\Big(\sum_{j=0}^\infty|d_j|z^j\Big)\Big(\sum_{j=0}^\infty z^j\Big),\\
&\hat{f}=\sum_{j=0}^\infty\big(|f_j|+|\theta_{-1}^2a_{j+2}|+|\theta_{-1}b_{j+1}|+|\theta_{-1}d_{j+2}|\big)z^j.
\end{align*}
These majorants are chosen so that the recurrence relation for $\hat{y}$ gives $|\theta_j|\leq\hat{\theta}_j$ for all $j$. For instance, the contribution of the term $d(z)y(z^2)$ to $P_n$ in (\ref{eq:nev:5}) is dominated by the contribution of $\hat{d}\hat{y}$. The same principle applies to a general $h$ and to (\ref{eq:nev:4}); in those cases we only give the analogous choices in (\ref{eq:appA:2}), without further explanation. The existence of $\hat{y}$ follows from \cite[p.300]{ahlfors1979complex}.

For a general $h$, analytic and non-zero at $z=0$, we set\footnote{The series $\sum_{j=0}^\infty(z\hat{h})^j$ converges for all $z$ in a small enough disc around $0$.}
\begin{align*}
&\hat{d}=\Big(\sum_{j=0}^\infty|d_j|z^j\Big)\Big(\sum_{j=0}^\infty(z\hat{h})^j\Big),\qquad \hat{h}=\sum_{j=0}^\infty|h_j|z^j,\\
&\hat{f}=\sum_{j=0}^\infty\big(|f_j|+|\theta_{-1}^2a_{j+2}|+|\theta_{-1}b_{j+1}|+|\theta_{-1}d_{j+2}|+|m_j|\big)z^j,
\end{align*}
where $\sum_{j=0}^\infty m_jz^j$ is the analytic part of $d(z)\theta_{-1}/(z^2h(z))$. A similar argument shows that $\hat{y}$ exists.

We now turn to (\ref{eq:nev:4}). First, assume that $\check{h}(z)\equiv q$ with $|q|\leq1$. Then (\ref{eq:nev:5}) can be written as $\Srec(n)\theta_{n+1}=P_n(\theta_{-1},\dots,\theta_{n})$, where
\begin{align*}
&\Srec(n)=n+1-2\theta_{-1}a_0-2\theta_{-1}d_0q^n-\theta_{-1}e_0q^{n+1}-\theta_{-1}e_0/q,\\
&P_{-1}=(a_1\theta_{-1}+b_0+c_0/q+d_1\theta_{-1}/q+e_1\theta_{-1}/q)\theta_{-1},
\end{align*}
and the $P_n$ are polynomials.\footnote{For a general $\check{h}$, the explicit form of $P_n$ is too lengthy to be useful, and we omit it.} Since $\Srec(n)\neq0$, there exists $K>0$ such that $|\Srec(n)|^{-1}\leq K$ for $n\geq0$. Let $\hat{y}=\sum_{j=0}^\infty\hat{\theta}_jz^j$ satisfy
\begin{equation}
\label{eq:appA:3}
\hat{y}=Kz(\hat{a}\hat{y}^2+\hat{b}\hat{y}+\hat{c}\hat{y}+\hat{d}\hat{y}^2+\hat{e}\hat{y}^2+\hat{f})+\hat{\theta}_0,
\end{equation}
where
\begin{align*}
&\hat{a}=\sum_{j=0}^\infty|a_j|z^j,\qquad \hat{b}=\sum_{j=0}^\infty\big(|b_j|+|2\theta_{-1}a_{j+1}|+|\theta_{-1}e_{j+1}|\big)z^j,\\
&\hat{c}=\sum_{j=0}^\infty\big(|c_j|+|2\theta_{-1}d_{j+1}/q|+|\theta_{-1}e_{j+1}/q|\big)z^j,\\
&\hat{d}=\sum_{j=0}^\infty|d_j|z^j,\qquad\hat{e}=\sum_{j=0}^\infty|e_j|z^j,\\
&\hat{f}=\sum_{j=0}^\infty\big(|f_j|+|\theta_{-1}^2a_{j+2}|+|\theta_{-1}b_{j+1}|+|\theta_{-1}c_{j+1}/q|\\
&\hspace{4.5em}+|\theta_{-1}^2d_{j+2}/q^2|+|\theta_{-1}^2e_{j+2}/q|\big)z^j.
\end{align*}
If $\hat{y}$ exists, then again $|\theta_j|\leq\hat{\theta}_j$ for $j\geq0$. Writing $\hat{y}=Y/z$, (\ref{eq:appA:3}) becomes
\begin{equation}
\label{eq:appA:4}
K(\hat{a}+\hat{d}+\hat{e})Y^2-Y+z\big(K(\hat{b}Y+\hat{c}Y+\hat{f})+\hat{\theta}_0\big)=0,
\end{equation}
and $Y$ is analytic at $0$ by \cite[p.300]{ahlfors1979complex}. Hence $\hat{y}$ exists.

For a general $\check{h}$, we modify the argument as follows. The recurrence relation (\ref{eq:nev:5}) takes the form
\[
\big(n+1-2\theta_{-1}a_0-2\theta_{-1}d_0\check{h}_0^n-\theta_{-1}e_0\check{h}_0^{n+1}-\theta_{-1}e_0/\check{h}_0\big)\theta_{n+1}=P_n(\theta_{-1},\theta_0,\dots,\theta_{n}).
\]
Let $\sum_{j=0}^\infty\tilde{e}_jz^j$ be the analytic part of $e_0\theta_{-1}\big(1/(z\check{h}(z))-1/(z\check{h}_0)\big)$, let $\sum_{j=0}^\infty\tilde{f}_jz^j=\theta_{-1}^2\sum_{j=0}^\infty e_{j+2}z^j/\check{h}(z)$, and let $\sum_{j=0}^\infty\check{f}_jz^j$ be the analytic part of
\[
\frac{e_0\theta_{-1}^2}{z^2\check{h}(z)}+\frac{e_1\theta_{-1}^2}{z\check{h}(z)}+\frac{c_0\theta_{-1}}{z\check{h}(z)}+\sum_{j=0}^\infty \frac{c_{j+1}\theta_{-1}z^j}{\check{h}(z)}.
\]
We set $\hat{\check{h}}=\sum_{j=0}^\infty|\check{h}_j|z^j$\footnote{Since $\check{h}$ is analytic at $0$ with $|\check{h}(0)|<1$, we have $\hat{\check{h}}(z)\leq\epsilon<1$ for all small enough $|z|$, so $\sum_{j=0}^\infty\hat{\check{h}}^j$ converges.} and
\begin{align*}
&\hat{b}=\sum_{j=0}^\infty\big(|b_j|+|2\theta_{-1}a_{j+1}|+|\tilde{e}_j|\big)z^j,\qquad \hat{c}=\Big(\sum_{j=0}^\infty|c_j|z^j\Big)\sum_{j=0}^\infty\hat{\check{h}}^j,\\
&\hat{d}=\Big|\frac{2d_0\theta_{-1}}{\check{h}_0}\Big|\Big(\sum_{j=0}^\infty|\check{h}_{j+1}|z^j\Big)\Big(\sum_{j=0}^\infty\hat{\check{h}}^j\Big)+\sum_{j=0}^\infty|\tilde{d}_j|z^j,
\end{align*}
where $\sum_{j=0}^\infty\tilde{d}_jz^j$ is the analytic part of $2d_0\theta_{-1}\big(1/(z\check{h}(z))-1/(z\check{h}_0)\big)$. Further, let
\begin{align*}
&\hat{e}_1=\Big(\sum_{j=0}^\infty|\theta_{-1}e_{j+1}|z^j\Big)\sum_{j=0}^\infty\hat{\check{h}}^j,\qquad \hat{e}_3=\Big(\sum_{j=0}^\infty|e_j|z^j\Big)\sum_{j=0}^\infty\hat{\check{h}}^j,\\
&\hat{e}_2=\sum_{j=0}^\infty|\check{e}_j|z^j,\quad\text{where }\sum_{j=0}^\infty\check{e}_jz^j=\sum_{j=0}^\infty \frac{e_{j+1}\theta_{-1}z^j}{\check{h}(z)},\\
&\hat{f}=\sum_{j=0}^\infty\big(|f_j|+|\theta_{-1}b_{j+1}|+|a_{j+2}\theta_{-1}^2|+|\tilde{f}_j|+|\check{f}_j|\big)z^j,
\end{align*}
and let $\hat{y}$ satisfy
\begin{equation}
\label{eq:appA:5}
\hat{y}=Kz\big(\hat{a}\hat{y}^2+\hat{b}\hat{y}+\hat{c}\hat{y}+\hat{d}\hat{y}^2+(\hat{e}_1+\hat{e}_2)\hat{y}+\hat{e}_3\hat{y}^2+\hat{f}\big)+\hat{\theta}_0.
\end{equation}
As before, $\hat{y}$ exists and $|\theta_j|\leq\hat{\theta}_j$ for all $j$.
\end{proof}

\begin{proof}[Proof of Corollary \ref{cor:nev:1}]
For (\ref{eq:nev:3}), $\Srec(n)=0$ for at most one $n$. For (\ref{eq:nev:4}), the set $\{n\in\mathbb{N}:\Srec(n)=0\}$ is finite, since $|\Srec(n)|\to\infty$ as $n\to\infty$; let $N$ be its largest element. To show that the analytic part of the Laurent series solution of (\ref{eq:nev:4}) converges, we only need to modify (\ref{eq:appA:5}) as follows; the proof for (\ref{eq:nev:3}) is similar. We again take $\hat{y}=\sum_{j=0}^\infty\hat{\theta}_jz^j$, now with $\hat{\theta}_j=|\theta_j|$ for $j=0,\dots,N$, and add terms to the right-hand side of (\ref{eq:appA:5}):
\begin{equation}
\label{eq:appA:6}
\hat{y}=Kz\big(\hat{a}\hat{y}^2+\hat{b}\hat{y}+\hat{c}\hat{y}+\hat{d}\hat{y}^2+(\hat{e}_1+\hat{e}_2)\hat{y}+\hat{e}_3\hat{y}^2+\hat{f}\big)+\hat{\theta}_0+\sum_{j=1}^N(\hat{\theta}_j-K_j)z^j,
\end{equation}
where $K_j$ is the coefficient of $z^j$ in the right-hand side of (\ref{eq:appA:5}). Then (\ref{eq:appA:6}) shows that $\hat{y}=\hat{\theta}_0+\hat{\theta}_1z+\cdots+\hat{\theta}_Nz^N+O(z^{N+1})$, with $\hat{\theta}_j=|\theta_j|$ for $j\leq N$ and $\hat{\theta}_j\geq|\theta_j|$ for $j>N$. The resonance condition at each positive integer root of $\Srec(n)$ need not hold automatically (see Theorem \ref{thm:qde:2}); it may impose conditions on $a,\dots,f$. In that case, the series exists only when these conditions hold.
\end{proof}

\section{An estimate for compositions of meromorphic functions}
\label{sec:appB}

The Valiron--Mohon'ko theorem \ref{thm:nev:7} implies that if $y$ is rational and $g$ is meromorphic, then $T(r,y\circ g)=(\deg y)\,T(r,g)+O(1)$. If instead $y$ is transcendental meromorphic, then $y\circ g$ is meromorphic whenever $g$ is entire. When $g$ is a polynomial of degree $\mu\geq2$ and $y$ is transcendental meromorphic, Jank and Volkmann \cite[§14]{jank1985einfuhrung} (see also \cite{rieppo2007class,fletcher2009meromorphic}) proved that there exist positive constants $d_1$ and $d_2$ such that
\begin{subequations}
\label{eq:appB:1}
\begin{align}
\label{eq:appB:2}
& N(r^\mu d_1,a,y)-O(1)\leq N(r,a,y\circ g)\leq N(r^\mu d_2,a,y)+O(1),\\
\label{eq:appB:3}
& (1-o(1))T(r^\mu d_1,y)\leq T(r,y\circ g)\leq (1+o(1))T(r^\mu d_2,y).
\end{align}
\end{subequations}
In fact the bounded $o(1)$ can be replaced by $O(1)$ as shown in \cite[Satz 14.1, 14.4]{jank1985einfuhrung} by choosing $d_1,d_2>0$ such that $d_1r^\mu\leq M(r,g)\leq d_2r^\mu$, we obtain
\begin{equation}
    \label{eq:appB:4}
    T(r^\mu d_1,y)-K_1\leq T(r,y\circ g)\leq T(r^\mu d_2,y)+K_2
\end{equation}
for some constants $K_1,K_2$. 

For the rest of this appendix, $g$ is a transcendental entire function, and
\[
M(r,g)=\max_{|z|=r}|g(z)|
\]
denotes its maximum modulus. P\'olya \cite{polya1926integral} proved that if $y$ is a transcendental entire function, then $M(r,y\circ g)\geq M(AM(r/2,g),y)$ for some constant $A>0$. This was later improved by Clunie \cite{clunie1955maximum} to
\begin{equation}
\label{eq:appB:5}
M(r,y\circ g)\geq M((1+o(1))M(r,g),y).
\end{equation}
A sharper bound, also due to Clunie \cite{clunie1970composition}, states that if $0<\alpha<1$ and $g(0)=0$, then $M(r,y\circ g)\geq M((1-\alpha)^2M(\alpha r,g)/(4\alpha),y)$.

In terms of the Nevanlinna characteristic, the corresponding estimates for a transcendental entire function $y$ can be found in \cite{clunie1970composition} (see also \cite[pp.51--54]{hayman1964meromorphic}):
\begin{equation}
\label{eq:appB:6}
\lim_{r\to\infty}\frac{T(r,y\circ g)}{T(r,y)}=\infty,\qquad \lim_{r\to\infty}\frac{T(r,y\circ g)}{T(r,g)}=\infty.
\end{equation}
In the same paper, Clunie also proved that, for a transcendental meromorphic function $y$,
\begin{equation}
\label{eq:appB:7}
\limsup_{r\to\infty}\frac{T(r,y\circ g)}{T(r,y)}=\infty,\qquad \lim_{r\to\infty}\frac{T(r,y\circ g)}{T(r,g)}=\infty.
\end{equation}
In fact, for any $K>0$, the set of $r$ with $T(r,y\circ g)\geq KT(r,y)$ has infinite measure.\footnote{This set contains a sequence of intervals $[a_n,b_n]$ with $b_n-a_n\to\infty$, so it can be combined with the LOTLD, whose exceptional set has finite measure.} The $\limsup$ in (\ref{eq:appB:7}) cannot be replaced by a limit in general, since Clunie constructed an example with $\liminf_{r\to\infty}T(r,y\circ g)/T(r,y)=0$; however, it can be replaced by a limit if $y$ or $g$ has finite order.

We now prove the following relation, which is used in the proof of Theorem \ref{thm:nev:9}:
\begin{equation}
\label{eq:appB:8}
\limsup\frac{T(r,y\circ g)}{T(r+\frac{1}{T(r,y\circ g)},y)}=\infty,
\end{equation}
The argument is a slight modification of Clunie's proof of \cite[Lemma 1, pp.82--85]{clunie1970composition}. Let $y$ be a transcendental meromorphic function and let $K>0$. There exist $R_0>0$ and an increasing sequence $R_n\to\infty$ such that, for $n\geq1$, all $r$ with $R_n\leq r\leq R_n^2$ and all $\omega$ with $R_0\leq|\omega|\leq r$, we have $n(r,1/(g-\omega))>K$. We may choose a value $a$ such that $N(r,a,y\circ g)=(1+o(1))T(r,y\circ g)$ and $N(r,a,y)=(1+o(1))T(r,y)$ (see \cite{ahlfors1930satz}, \cite[p.276]{nevanlinna1970analytic} and \cite[p.9]{hayman1964meromorphic}), and we may assume that $a=0$. Then, for $R_n^2/2+1\leq r\leq 9R_n^2/10$,
\begin{equation*}
    n(r,1/(y\circ g))\geq K\big(n(r,1/y)-n(R_0,1/y)\big),
\end{equation*}
and hence
\begin{align*}
N(r,1/(y\circ g))&\geq K\int_{R_n}^r\frac{n(t,1/y)}{t}\,dt+O(\ln R_n)\\
&\geq K\big(N(r,1/y)-N(R_n,1/y)\big)+O(\ln R_n).
\end{align*}
Thus $T(r,y\circ g)\geq (K+o(1))\big(T(r,y)-T(R_n,y)\big)$.

Next, we take $r$ in $[R_n^2/2+1,\,9R_n^2/10]$ of the form $r=R_n'+1/T(R_n',y\circ g)$, outside the exceptional set in Borel's lemma \cite[Lemmas 1.1 and 2.1]{halburd2010nondecreasing}: for any positive non-decreasing function $F$,
\begin{equation}
\label{eq:appB:9}
F(r+1/F(r))<2F(r)
\end{equation}
outside a set of finite linear measure. Such a choice is possible, since $9R_n^2/10-R_n^2/2-1\geq R_n^2/3$ for all large $n$. Applying (\ref{eq:appB:9}) with $F(r)=T(r,y\circ g)$, we obtain
\begin{equation*}
    \begin{aligned}T\Big(R_n'+\frac{1}{T(R_n',y\circ g)},y\circ g\Big)\geq{}& (K+o(1))\Big(T\Big(R_n'+\frac{1}{T(R_n',y\circ g)},y\Big)-T(R_n,y)\Big).\end{aligned}
\end{equation*}
By (\ref{eq:appB:9}), $T(R_n'+\frac{1}{T(R_n',y\circ g)},y\circ g)<2T(R_n',y\circ g)$, and thus
\begin{equation*}
    T(R_n',y\circ g)\geq\frac{K}{2}\left(T\left(R_n'+\frac{1}{T(R_n',y\circ g)},y\right)-T(R_n,y)\right).
\end{equation*}
By Theorem \ref{thm:nev:4}, $(2+o(1))T(R_n,y)\leq T(R_n^2/2,y)\leq T(R_n'+\frac{1}{T(R_n',y\circ g)},y)$, which yields
\begin{equation*}
    T(R_n',y\circ g)\geq\frac{K}{5}T\left(R_n'+\frac{1}{T(R_n',y\circ g)},y\right)
\end{equation*}
Hence $\liminf_{n\to\infty} T(R_n',y\circ g)/T(R_n'+1/T(R_n',y\circ g),y)\geq K/5$. Since $K$ is arbitrary, this proves (\ref{eq:appB:8}).

\bibliographystyle{plain}
\bibliography{bibliography3}

\end{document}